\documentclass[number,citesort,seceqn,dvips]{arxbj}
\usepackage{mathrsfs}
\aid{0}
\volume{17}
\issue{3}
\pubyear{2011}
\firstpage{987}
\lastpage{1014}
\doi{10.3150/10-BEJ307}

\makeatletter
\newtheorem{theorem}{Theorem}[section]
\newtheorem{proposition}[theorem]{Proposition}
\newtheorem{lemma}[theorem]{Lemma}

\newcommand{\eqref}[1]{(\ref{#1})}

\newtheorem{geo}{Proposition}[section]

\newtheorem{qiterates}[geo]{Proposition}
\newtheorem{kiterates}[geo]{Lemma}
\newtheorem{sllnu}[geo]{Lemma}
\newtheorem{coupu}[geo]{Proposition}

\makeatother

\begin{document}
\begin{frontmatter}

\title{On nonlinear Markov chain Monte Carlo}
\runtitle{On nonlinear Markov chain Monte Carlo}

\begin{aug}
\author[1]{\fnms{Christophe} \snm{Andrieu}\thanksref{1}\ead[label=e1]{c.andrieu@bris.ac.uk}},
\author[2]{\fnms{Ajay} \snm{Jasra}\thanksref{2}\ead[label=e2]{a.jasra@ic.ac.uk}},
\author[3]{\fnms{Arnaud} \snm{Doucet}\thanksref{3}\ead[label=e3]{a.doucet@cs.ubc.ca}}
\and
\author[4]{\fnms{Pierre} \snm{Del Moral}\thanksref{4}\ead[label=e4]{pierre.del-moral@inria.fr}}

\runauthor{Andrieu, Jasra, Doucet and Del Moral}
\address[1]{Department of Mathematics, University of Bristol, Bristol
BS8 1TW, UK.\\ \printead{e1}}
\address[2]{Department of Mathematics, Imperial College London,
London, SW7 2AZ, UK.\\ \printead{e2}}
\address[3]{Department of Statistics, University of British Columbia,
Vancouver, V6T 1Z2, Canada.\\ \printead{e3}}
\address[4]{Centre INRIA \& Institut de Math\'ematiques de Bordeaux,
Universit\'e de Bordeaux I, 33405, France.\\ \printead{e4}}
\end{aug}

\received{\smonth{4} \syear{2007}}
\revised{\smonth{3} \syear{2010}}

%
\begin{abstract}
Let $\mathscr{P}(E)$ be the space of probability measures on a
measurable space $(E,\mathcal{E})$. In this paper we introduce a class
of nonlinear Markov chain Monte Carlo (MCMC) methods for simulating
from a probability measure $\pi\in\mathscr{P}(E)$. Nonlinear Markov
kernels (see [\textit{Feynman--Kac Formulae: Genealogical and
Interacting Particle Systems with Applications} (2004) Springer])
$K\dvtx\mathscr{P}(E)\times E\rightarrow\mathscr{P}(E)$ can be
constructed to, in some sense, improve over MCMC methods. However, such
nonlinear kernels cannot be simulated exactly, so approximations of the
nonlinear kernels are constructed using auxiliary or potentially
self-interacting chains. Several nonlinear kernels are presented and it
is demonstrated that, under some conditions, the associated
approximations exhibit a strong law of large numbers; our proof
technique is via the Poisson equation and Foster--Lyapunov conditions.
We investigate the performance of our approximations with some
simulations.
\end{abstract}

%

\end{frontmatter}

\section{Introduction}\label{1}

Monte Carlo simulation is one of the most important elements of
computational statistics. This is because of its relative simplicity
and computational convenience in constructing estimates of
high-dimensional integrals. That is, for a $\pi$-integrable
$f\dvtx E\rightarrow\mathbb{R}$, we approximate:
%
\begin{equation}\label{integ}
\pi(f)  :=  \int_{E}f(x)\pi(\mathrm{d}x)
\end{equation}
by
\[
S_n^{X}(f)  =  \frac{1}{n+1}\sum_{i=0}^{n}f(X_i),
\]
where $S_n^X(\mathrm{d}u):=\frac{1}{n+1}\sum_{i=0}^n\delta_{X_i}(\mathrm{d}u)$ is the
empirical measure based upon random variables $\{X_k\}_{0\leq k\leq n}$
drawn from $\pi$. Such integrals appear routinely in Bayesian
statistics, in terms of posterior expectations; see \cite{robert} and
the references therein. In those cases, $E$ is often of very high
dimension and complex simulation methods such as MCMC \cite{robert}
and sequential Monte Carlo (SMC) \cite{delmoral1,doucet} need to be
used.

It has long been known by Monte Carlo specialists that standard MCMC
algorithms often have difficulties in simulating from complicated
distributions -- for example, when the target $\pi$ exhibits multiple
modes and/or possesses strong dependencies between subcomponents of
$X$. In the former case, the Markov chain can take an unreasonable
amount of time to jump between these modes and the estimates of
(\ref{integ}) are very inaccurate.

As a result, there have been a large number of alternative methods
proposed in the literature; we detail some of them here. Many of these
approaches have relied upon MCMC techniques such as adaptive MCMC
\cite{andrieu1,haario}, which, in some instances, attempts to improve
the mixing properties of the transition kernel by using the information
learned in the past. In addition, there are methods that rely upon the
simulation of parallel Markov chains \cite{geyer} and genetic algorithm
type moves; see \cite{jasra} for a review. These latter methods use the
idea of running some of the parallel chains with invariant probability
measure $\eta$, where $\eta$ is easier to explore and is related to
$\pi$; hence the samples of the parallel chains can provide valuable
information for simulating from $\pi$. Extensions to MCMC-based
simulation methods have combined MCMC with SMC ideas, see, for example,
\cite{pmcmc,delmoral3}. Such approaches are often more flexible than
MCMC.

In this paper, we consider another alternative: nonlinear MCMC via
auxiliary or
self-interacting approximations. Such methods rely primarily upon the
ideas of MCMC. However, it is demonstrated below that the auxiliary/self-interacting approximation idea is similar to that of approximating
Feynman--Kac formulae \cite{delmoral1} and as such is linked to SMC
methodology. It should be noted that related ideas
have appeared, directly in \cite{brockwell} and indirectly in
\cite{kou}; see \cite{andrieu5,atchade2} for some theoretical
analysis. Subsequent to the first versions of this work
\cite{nlfirst} a variety of related articles have appeared:
\cite{atchade1,atchade2,atchadecheatingafuck}; we cite these where
appropriate, but note the substantial overlap between our work and
these papers.

\subsection{Nonlinear Markov kernels via interacting approximations}\label{11}

Standard MCMC algorithms rely on Markov kernels of the form
$K\dvtx E\rightarrow\mathscr{P}(E)$. These Markov kernels are \textit{linear}
operators on $\mathscr{P}(E)$; that is, ~$\mu(\mathrm{d}y)=\int_E
\xi(\mathrm{d}x)K(x,\mathrm{d}y),$ where $\mu,\xi\in\mathscr{P}(E)$. A \textit{nonlinear}
Markov kernel $K\dvtx \mathscr{P}(E)\times E\rightarrow\mathscr{P}(E)$ is
defined as a nonlinear operator on the space of probability measures.
Nonlinear Markov kernels, $K_\mu$, can often be constructed to exhibit
superior mixing properties to ordinary MCMC versions. For example, let
%
\begin{equation} \label{eq:nonlin1}
K_{\mu}(x,\mathrm{d}y)  =  (1-\epsilon)K(x,\mathrm{d}y) +
\epsilon\int_{E}\mu(\mathrm{d}z)K(z,\mathrm{d}y),
\end{equation}
where $K$ is a Markov kernel of invariant distribution $\pi$,
$\epsilon\in(0,1)$ and $\mu\in\mathscr{P}(E)$. Simulating from
$K_\pi$
is clearly desirable as we allow regenerations from $\pi$, with $K_\pi$
strongly uniformly ergodic (see~\cite{roberts}). However, in most
cases, it is not possible to simulate from $K_\pi$ and, instead, an
approximation is proposed.

A self-interacting Markov chain (see~\cite{delmoral4}) generates a
stochastic process $\{X_n\}_{n\geq0}$ that is allowed to interact with
values realized in the past. That is, we might approximate, at time
$n+1$, $\mu$ by $S_n^X$. This process corresponds to generating a value
from the history of the process, and then a mutation step, via the
kernel $K$. In practice, the self-interaction can lead to very poor
algorithmic performance \cite{nlfirst}; an auxiliary Markov chain is
used to approximate the nonlinear kernel.

\subsection{Motivation and structure of the paper}\label{12}

In the context of stochastic simulation, self-interacting Markov chains
(SIMCs), or IMCs, can be thought of as storing modes and then allowing
the algorithm to return to them in a relatively simple way. Parametric
adaptive MCMC can be thought of as an indirect application of this
idea, where parameters of the kernel are optimized via a stochastic
approximation algorithm. This approach does not retain all of the
features of previously visited states. In other words, SIMCs can be
considered as a nonparametric, or infinite-dimensional, generalization
of parametric adaptive MCMC. It is thus the attractive idea of being
able to fully exploit the information provided by the previous samples
that has motivated us to investigate such algorithms.

This paper is structured as follows. We begin by giving our notation in
Section \ref{2}. In Section~\ref{3} our simulation methods are described and
several nonlinear Markov kernels and self-interacting approximations
are introduced. In Section \ref{4} we introduce some assumptions and some
preliminary results, which are used to prove a strong law of large
numbers (SLLN). In Sections \ref{5} and~\ref{6}, some technical proofs and the SLLN
are presented; this is for a particular nonlinear kernel introduced in
Section \ref{3}. This analysis is of interest from a theoretical point of
view: it brings together the literature of measure-valued processes and
interacting particle systems \cite{delmoral1} used in SMC and the
relatively recent literature on general state space Markov chains
\cite{meyn} used in MCMC. In Section \ref{7} some algorithms are
investigated; our assumptions are verified and some parameter settings
are investigated for a toy example. In Section \ref{sec:summary} some extensions to our
ideas are discussed. The proofs are all given in the \hyperref[a1]{Appendices}.

\section{Notation and definitions}\label{2}

\subsection{Notation}\label{21}

\subsubsection{Probability and measure}\label{211}
Define a measurable space $(E,\mathcal{E})$. Throughout, $\mathcal{E}$
will be assumed countably generated. $\mathscr{B}(\mathbb{R}^k)$,
$k\in\mathbb{N}$ is used to represent the Borel sets with Lebesgue
measure denoted by $\mathrm{d}x$.

For a stochastic process $\{X_n\}_{n\geq0}$ on
$(E^{\mathbb{N}},\mathcal{E}^{\otimes\mathbb{N}})$,
$\mathcal{G}^X_n=\sigma(X_0,\ldots,X_n)$ denotes the natural filtration.
$\mathbb{P}_{\mu}$ is taken as a probability law of a stochastic
process with initial distribution $\mu$ and $\mathbb{E}_{\mu}$ the
associated expectation. If $\mu=\delta_x$, with $\delta$ the Dirac
measure, $\mathbb{P}_x$ (resp.,~$\mathbb{E}_{x}$) is used instead
of $\mathbb{P}_{\delta_x}$ (resp.,~$\mathbb{E}_{\delta_x}$). For
$\mu\in\mathscr{P}(E)$, the product measure is written
$\mu\times\mu=\mu^{\otimes2}$, with a clear generalization to higher
order products. For measurable $f\dvtx E\rightarrow\mathbb{R}$,
$\mu(f)=\int_E f(x)\mu(\mathrm{d}x)$.

If a $\sigma$-finite measure $\pi$ is dominated by another $\eta$
(denoted $\pi\ll\eta$), the Radon--Nikodym derivative is written with
the same notation (e.g.,~if $\pi\ll\eta$, then
$\pi(x)/\eta(x)=\mathrm{d}\pi/\mathrm{d}\eta(x)$). For $\sigma$-finite measures $\pi
$ and
$\eta$, $\pi\sim\eta$ denotes mutual absolute continuity.

\subsubsection{Markov chains}\label{212}

Let $(E,\mathcal{E})$ be a measurable space. Throughout for a Markov
transition kernel $K\dvtx E\rightarrow\mathscr{P}(E)$ the following standard
notation is used: for measurable $f\dvtx E\rightarrow\mathbb{R}$,
$K(f)(x):=\int_{E}f(y)K(x,\mathrm{d}y)$ and for $\mu\in\mathscr{P}(E)$ $\mu
K(f):=\int_{E}K(f)(x)\mu(\mathrm{d}x)$.

For $K_{\mu}$, $K\dvtx E\times\mathscr{P}(E) \rightarrow\mathscr{P}(E)$,
given its existence, we will denote by $\omega(\mu)$
($\omega\dvtx \mathscr{P}(E)\rightarrow\mathscr{P}(E)$) the invariant
distribution of this Markov kernel. Recall that the empirical measure
of an arbitrary stochastic process
$(E^{\mathbb{N}},\mathcal{E}^{\mathbb{N}},\{X_n\}_{n\geq
0},\mathbb{P})$ is defined, at time $n$, as
%
\begin{equation} \label{eq:def_empirical_measure}
S_n^{X}(\mathrm{d}u) := \frac{1}{n+1}\sum_{i=0}^n\delta_{X_i}(\mathrm{d}u).
\end{equation}

Throughout this paper, we are concerned with two nonlinear kernels of
the form
\begin{eqnarray*}
K_{\mu}(x,\mathrm{d}y)
& = &
(1-\epsilon) K(x,\mathrm{d}y) + \epsilon\Phi(\mu)(\mathrm{d}y),
\\
\Phi(\mu)(f)
& = &
\int_E \frac{g(y)f(y)}{\mu(g)}\mu(\mathrm{d}y),
\end{eqnarray*}
where $K\dvtx E\rightarrow\mathscr{P}(E)$, $F\dvtx E\times\mathscr{P}(E)
\rightarrow\mathscr{P}(E)$ (see \cite{delmoral1} for more on $\Phi$)
and
\begin{eqnarray}\label{eq:generic_K_mu}
K_{\mu}(x,\mathrm{d}y)
& = &
(1-\epsilon) K(x,\mathrm{d}y) + \epsilon Q_{\mu}(x,\mathrm{d}y),\nonumber
\\[-8pt]\\[-8pt]
Q_{\mu}(f)(x)
& = &
\int_{E} \mu(\mathrm{d}u) \alpha(x,u)[f(u)-f(x)] +f(x),\nonumber
\end{eqnarray}
where $\alpha(x,u)$ is defined later on.

\subsubsection{Norms}\label{213}

For any $k\in\mathbb{N}$, the Euclidean norm of $x \in\mathbb{R}^{k}$
is denoted $|x|$. For $f\dvtx E\rightarrow\mathbb{R}^{n}$, $n\in
\mathbb{N}$, $|f|_{\infty}:=\sup_{x\in E}|f(x)|$. For
$f:E\rightarrow\mathbb{R}^{n}$ the $\mathbb{L}_p$-norm is defined,
assuming it exists, as $(\int_E |f(x)|^{p}\,\mathrm{d}\mu)^{1/p}$ for
$\mu\in\mathscr{P}(E)$. For $V\dvtx E\rightarrow[1,\infty)$ and
$f\dvtx E\rightarrow\mathbb{R}^{n}$
\[
|f|_V  := \sup_{x\in E}\frac{|f(x)|}{V(x)} .
\]
$\mathscr{L}_V$ is the class of functions $f\dvtx E\rightarrow
\mathbb{R}^{n}$ such that $|f|_V<\infty$. We also use the notions of
the $V$-total variation for a signed measure
\[
\|\lambda\|_V := \sup_{|f|\leq V}|\lambda(f)| ,
\]
and the $V$-norm operator between two kernels
$K_1,K_2\dvtx E\rightarrow\mathscr{P}(E)$
\[
\Vert\!\vert K_1-K_2\vert\!\Vert_V  :=  \sup_{x\in
E}\frac{\|K_1(x,\cdot)-K_2(x,\cdot)\|_V}{V(x)}.
\]

\subsubsection{Miscellaneous}\label{214}

The notation $a\vee b:=\max\{a,b\}$ (resp.,~$a\wedge
b:=\min\{a,b\}$) is adopted. The indicator function of $A\subset E$ is
written $\mathbb{I}_A(x)$. $\mathbb{N}_0=\mathbb{N}\cup\{0\}$.
Throughout the paper we denote~a~gene\-ric finite constant as $M$, that
is, the value of $M$ may change from line to line in the proofs and is
local to each proof.

\section{Nonlinear MCMC}\label{3}

\subsection{Nonlinear Markov kernels}
\label{nlmcmc}

Nonlinear MCMC can be characterised by the following procedure:
\begin{itemize}
\item Identify a nonlinear kernel that admits $\pi$ as an invariant
distribution and can be expected
to mix faster than an ordinary MCMC kernel; for example,~(\ref
{eq:nonlin1}).
\item Construct a stochastic process that approximates the kernel,
which can be simulated in practice.
\end{itemize}

Based upon the previous work \cite{nlfirst}, we consider auxiliary
stochastic processes to approximate the nonlinear kernel. That is, it
has been found in \cite{nlfirst} that using the past history to
approximate the nonlinear kernel leads to very poor performance. All of
the processes that are simulated in this paper use an auxiliary Markov
chain to approximate the nonlinear kernel. The difficulty is then to
design sensible nonlinear kernels that may lead to good empirical
performance. The two kernels we have designed are below.

\subsection{Selection/mutation with potential}\label{sec:selection_mut}
Let $P$ be an MCMC kernel of invariant distribution $\eta$, and assume
$\pi\ll\eta$. Let $g(v) = \frac{\pi(v)}{\eta(v)}$ and set $K$ to
be an
MCMC kernel of invariant distribution $\pi$. Consider the nonlinear
kernel
\[
K_{\mu}(x,\mathrm{d}x') = (1-\epsilon) K(x,\mathrm{d}x') + \epsilon\Phi(\mu)(\mathrm{d}x');
\]
clearly, if $\mu=\eta$, then one has $\pi K_{\eta} = \pi$.

If it is possible to sample exactly from $\eta$, then one could sample
exactly from $K_{\eta}$. However, for efficient algorithms, this will
not be the case. The following approximation is adopted at time-step
$n+1$ of the simulation:
\[
[(1-\epsilon)K(x_{n},\mathrm{d}x_{n+1}) +
\epsilon\Phi(S_n^Y)(\mathrm{d}x_{n+1})]P(y_{n},\mathrm{d}y_{n+1});
\]
that is, we are `feeding' the chain $\{X_n\}_{n\geq0}$ the empirical
measure $S_n^Y$. Intuitively, as $n$ grows large,
$S_n^Y(f)\rightarrow\eta(f)$ and one samples from the original kernel
of interest.

\subsection{Auxiliary self-interaction with genetic moves}\label{sec:auxselfinteration}
For any $\mu\in\mathscr{P}(E)$ we
define a nonlinear Markov kernel $Q_\mu\dvtx \mathscr{P}(E)\times
E\rightarrow\mathscr{P}(E)$
\[
Q_{\mu}(f)(x) = \int_{E} \mu(\mathrm{d}u) \alpha(x,u)[f(u)-f(x)] + f(x)
\]
and for $\pi\sim\eta$
\[
\alpha(x,y) = 1 \wedge\frac{\pi(y)\eta(x)}{\pi(x)\eta(y)}.
\]
The idea here is to generate a sample from $\mu$ and accept or reject
it as the new state on the basis of the probability $\alpha$. Clearly,
$\pi Q_{\eta} = \pi$. Letting $K$ and $P$ be as above, the process is
simulated according to
\[
\{(1-\epsilon)K(x_n,\mathrm{d}x_{n+1})+\epsilon
Q_{S_n^Y}(x_n,\mathrm{d}x_{n+1})\}P(y_n,\mathrm{d}y_{n+1})
\]
at time $n+1$.

\subsection{Some comments}\label{34}

In the example in Section \ref{sec:selection_mut} we attempt to use
some measure of information, through~$g$, to assist the resampling. The
example of Section \ref{sec:auxselfinteration} provides a way to
control the information that is provided by the approximation $S_n^Y$.
That is, the kernel $Q_{S_n^Y}$, via $\alpha$ and the possible
rejection, will provide a criterion to check the consistency with the
target of the value drawn from $S_n^Y$. This may help improve
estimation, if $S_n^Y$ converges slowly. Note that the algorithm is
related to, but less sophisticated than, that of \cite{kou}. This is
because we do not consider exchanges to occur between states in
equi-energy rings.

It should be remarked that similar kernels are investigated in
\cite{atchade2}. The author deduces that for a toy example it is hard
to justify the use of such adaptive methods. However, a~potential
criticism of that study is that it is for a unimodal target;
`advanced' methods are seldom necessary for such scenarios. This is
discussed further in Section \ref{sec:toyexample}.

\subsection{Algorithm}\label{sec:algo_descrip}

The algorithm is (with the appropriate $\Phi(\mu)$ or $Q_{\mu}$):
\begin{enumerate}[0.]
\item[0.] (Initialization): Set $n=0$ and $X_0=x$, $Y_0=y$,
$S^{Y}_{0}=\delta_y$.
\item[1.] (Iteration): Set $n=n+1$, simulate $Y_n\sim
P(Y_{n-1},\cdot)$ and $X_n\sim
K_{S_{n-1}^Y}(X_{n-1},\cdot)$.
\item[2.] (Update): $S_n^Y = S_{n-1}^Y + \frac{1}{n+1}[\delta_{Y_n}
- S_{n-1}^Y]$ and return to 1.
\end{enumerate}

\section{Assumptions}\label{4}


We now seek to prove an SLLN for the nonlinear MCMC algorithm described
in Section~\ref{sec:auxselfinteration}. Recall that we simulate a
stochastic process on $((E\times E)^{\mathbb{N}},(\mathcal{E}\otimes
\mathcal{E})^{\otimes\mathbb{N}}, \{X_n,Y_n\}_{n\geq0},$
$\{\mathcal{G}_n\}_{n\geq0},\mathbb{P}_{(x,y)})$, $(x,y)\in E\times
E$,
with
finite-dimensional law:
\[
\mathbb{P}_{(x,y),n}(d(x_0,y_0,\dots,x_n,y_n))  =
\delta_{(x,y)}(d(x_0,y_0))
\prod_{i=0}^{n-1}K_{S_i^y}(x_i,\mathrm{d}x_{i+1})P(y_i,\mathrm{d}y_{i+1})
.
\]
Note that the natural filtration is denoted as
$\mathcal{G}_n=\mathcal{G}_n^{X,Y}$ for notational simplicity. Since
$\{Y_n\}$ is generated independently of $\{X_n\}$, we denote the
probability law of the Markov chain $\{Y_n\}$ as $\mathbb{Q}_y$. Note,
again, that the proofs are given in the \hyperref[a1]{Appendices}.

\subsection{Assumptions}\label{sec:assumptions}
Our assumptions on $K$, used to define our
process, are now given. For $\bar{M}\in\mathbb{R}_+$, the notation
$\mathscr{P}_{\bar{M}}(E)=\{\mu\in\mathscr{P}(E)\dvtx \mu(V)<\bar{M}\}
$ is
adopted, with $V$ defined below. In the remainder of the paper we say
that a set $C \subset E$ is $(1,\theta)$-small if it satisfies a 1-step
minorization condition, with parameter $\theta\in(0,1)$.
\begin{enumerate}[(A1)]
\item[(A1)] Stability of $K$.
\begin{enumerate}[(iii)]
\item[(i)](\textit{Invariance and irreducibility}).
$K\dvtx E \rightarrow\mathscr{P}(E)$ is a $\pi$-invariant and~$\phi$-irre\-ducible Markov kernel.
\item[(ii)](\textit{One-step minorization on
level sets}). Define $C_d:=\{x\in
E\dvtx V(x)\leq d\}$ for any $d\in(1,\infty)$. We assume that for any
$d\geq
1$, $C_d$ is $(1,\theta_{d})$-small for some $\theta_d\in(0,1)$ and
$\nu_{d}\in\mathscr{P}(E)$.
\item[(iii)](\textit{One-step drift condition}). There exist
$V\dvtx E\rightarrow[1,\infty)$ such that\break $\lim_{|x|\rightarrow\infty}V(x)
= \infty$, $\lambda<1$, $b<\infty$, $C \in\mathcal{E}$ such that for
any $x \in E$
\[
KV(x) \leq \lambda V(x) + b\mathbb{I}_{C}(x) .
\]
%
\end{enumerate}
%

\item[(A2)]\textit{Stability of $P$.}
\begin{enumerate}[(iii)]
\item[(i)]($W$-\textit{uniform ergodicity}).  $P\dvtx E
\rightarrow\mathscr{P}(E)$ is an $\eta$-invariant Markov kernel.
Furthermore, there exists $W\dvtx E\rightarrow[1, \infty)$ such that $P$ is
a $W$-uniformly ergodic Markov transition kernel with a one-step drift
condition and one-step minorization condition. In addition, there
exists an $r^*\in(0,1]$ such that $V\in\mathscr{L}_{W^{r^*}}$ (where
$V\dvtx E \rightarrow[1,\infty)$ is defined in
(A1)(iii)).
\end{enumerate}
\item[(A3)]\textit{State-space constraint}
\[
(E,\mathcal{E}) \mbox{ is Polish}.
\]
\end{enumerate}

\subsection{Discussion of the assumptions}

Our proofs of the SLLN will rely upon a martingale approximation via
the solution of the Poisson equation (e.g.,~\cite{glynn}). For any
$\bar{M}<\infty$, (A1) will allow us to
establish a drift condition for the kernel $K_{\mu}$ that is uniform in
$\mu\in\mathscr{P}_{\bar{M}}(E);$ see \cite{andrieu1}. In turn, one
can establish: the existence of a solution to Poisson's equation, the
existence of an invariant measure $\omega(\mu)$ for $K_{\mu}$ and
regularity properties uniform in $\mu\in\mathscr{P}_{\bar{M}}(E)$.
Then, due to~(A2), the following facts are
exploited: $\{S_n^Y(V)\}$ is $\mathbb{Q}_y$-a.s. finite and given
$\{S_n^Y(V)\}$, $\{X_n\}$ is a~Markov chain.
(A1) and (A.2) appear quite
strong, but can be verified in some important cases such as for random
walk Metropolis kernels; see \cite{jarner}, for example.

A key result, relying on both (A2) and
(A3), which is of interest in itself, is that of the
$\mathbb{Q}_y$-a.s. convergence of $V$-statistics of $\{Y_i\}$. This
result will enable us to show that, $\mathbb{Q}_y$-a.s.,
$\omega(S_i^Y)\rightarrow\omega(\eta)$; this is needed for our proof.


\section{Common properties of $K_{\mu}$}\label{5}

Using standard drift and minorization conditions, the existence of an
invariant probability measure is established for any
$\mu\in\mathscr{P}_{\infty}(E)$ under (A1).

\begin{proposition}\label{prop:K_mu_is_geometric}
Assume (\textup{A1}). Let $\epsilon\in(0,1)$ as in
\eqref{eq:generic_K_mu}, $\bar{M}\in(0,\infty)$, then for $d>
\epsilon
\bar{M} /[(1-\epsilon)(1-\lambda)]$ with $\lambda$ and $b$ as in
\textup{(A1)(iii)}:
\begin{enumerate}
\item There exist $(\theta_{d}',\nu_{d}) \in(0,1)\times\mathscr{P}(E)$
such that for any $\mu\in\mathscr{P}_{\bar{M}}(E)$ and $(x,A)\in
E\times\mathcal{E}$:
\begin{eqnarray*}
K_\mu(x,A) & \geq& \mathbb{I}_{C_{d}}(x)\theta_{d}'\nu_{d}(A),
\\
K_\mu V(x) & \leq& \tilde{\lambda}V(x) + \tilde{b}
\mathbb{I}_{C_{d}}(x)
\end{eqnarray*}
%
with $\tilde{\lambda}=(1-\epsilon)\lambda+ \epsilon+ \frac
{\epsilon
\bar{M}}{d}<1$, $\tilde{b}=(1-\epsilon)[\lambda d + b] + \epsilon
[\bar{M} + d]$.
\item There exists a function
$\omega\dvtx \mathscr{P}_{\infty}(E)\rightarrow\mathscr{P}_{\infty}(E)$,
such that for any $\mu\in\mathscr{P}_{\infty}(E)$
\[
\omega(\mu) = \omega(\mu)K_{\mu}.
\]
%
%
\item There exist constants, $\rho\in(0,1)$ and $M<\infty$ depending
upon $\bar{M}$, $\epsilon$, $\lambda$, $b$, $V$, d, $\theta_d$ (as
defined in equation (\ref{eq:generic_K_mu}) and
\textup{(A1)}), such that for any
$\mu\in\mathscr{P}_{\bar{M}}(E)$, $r\in(0,1]$ and
$f\in\mathscr{L}_{V^r}$
\[
|K_{\mu}^n(f) - \omega(\mu)(f)|_{V^r}  \leq M|f|_{V^r}\rho^n .
\]
\end{enumerate}
\end{proposition}

Some continuity properties associated with the invariant measures are
as follows.

\begin{proposition} \label{prop:invmeasnl3}
Assume (\textup{A1}) and let $\bar{M}\in(0,\infty)$.
Then there exists $M<\infty$ (depending solely on $\bar{M}$ and the
constants in (\textup{A1})) such that for any
$r\in(0,1]$, $\mu,\xi\in\mathscr{P}_{\bar{M}}(E)$,
\begin{eqnarray*}
\|\omega(\xi)-\omega(\mu)\|_{V^r} & \leq& M\Vert\!\vert K_{\xi} -
K_{\mu}\vert\!\Vert_{V^r} .
\end{eqnarray*}
\end{proposition}

Noting that for any $\mu,\xi\in\mathscr{P}(E)$ and $r\in[0,1]$,
$\Vert\!\vert K_{\xi} - K_{\mu}\vert\!\Vert_{V^r}=\epsilon\Vert\!\vert Q_\xi- Q_\mu\vert\!\Vert_{V^r}$ we
establish global Lipschitz continuity results for $\mu\mapsto Q_\mu$,
which, together with the result above, will allow us to deduce uniform
Lipschitz continuity of $\mu\rightarrow K_{\mu}$ on
$\mathscr{P}_{\bar{M}}(E)$ for any $\bar{M}\in(0,\infty)$. This is to
be used in the proofs of many of the subsequent results.

\begin{proposition} \label{prop:lipschitznl3}
Let $\mu, \xi\in\mathscr{P}_{\infty}(E)$, then for any $r\in(0,1]$:
\[
\Vert\!\vert Q_{\mu} - Q_{\xi}\vert\!\Vert_{V^r} \leq2\|\mu-\xi\|_{V^r}.
\]
\end{proposition}

\section{Law of large numbers}\label{6}

\subsection{Main result}\label{61}

Our main result is the following SLLN.

\begin{theorem} \label{llnnl3}
Assume \textup{(A1)--(A3)}. Let $r\in[0, 1)$. Then for any
$f\in\mathscr{L}_{V^r}$, $(x,y)\in E\times E$
\[
S_n^X(f)\stackrel{a.s.}{\longrightarrow}_{\mathbb{P}_{(x,y)}} \pi(f)
.
\]
\end{theorem}

The proof is detailed in Appendix \ref{app:proofmainresult}, but we
outline its main steps below.

\subsection{Strategy of the proof}\label{sec:strategy}

The strategy of the proof is now outlined. Introduce the following
sequence of probability distributions
$\{S^{\omega}_n:=1/(n+1)\sum_{i=0}^n\omega(S_i^Y)\}_{n\geq0}$, where
$\omega(\mu)$ is the invariant measure of~$K_\mu$ (which, if
$\mu=S_m^Y$, exists $\mathbb{Q}_y$-a.s.). This distribution can be used
as a re-centering term in the following decomposition,
%
\begin{equation} \label{eq:prfdecomp}
S_n^X(f) -\pi(f) = S_n^X(f) - S^{\omega}_n(f) + S^{\omega}_n(f)
-\pi(f) .
\end{equation}
Let $\mu\in\{S_n^Y(f)\}$ and assume, for now, the almost sure
existence of a solution $\hat{f}_{\mu}$ to Poisson's equation, that is,
such that for any $x\in E$
\[
f(x)-\omega(\mu)(f)  =  \hat{f}_{\mu}(x) - K_{\mu}(\hat{f}_{\mu})(x).
\]
Then, the first term on the right-hand side of (\ref{eq:prfdecomp}) can
be rewritten as
\begin{eqnarray} \label{eq:decomp}
(n+1)[S_n^X-S_n^{\omega}](f)  &=&  M_{n+1} +
\sum_{m=0}^{n}[\hat{f}_{S_{m+1}^Y}(X_{m+1})-\hat{f}_{S_m^Y}(X_{m+1})]\nonumber
\\[-8pt]\\[-8pt]
&&{}+
\hat{f}_{S_0^Y}(X_0) - \hat{f}_{S_{n+1}^Y}(X_{n+1}) ,\nonumber
\end{eqnarray}
where
\[
M_{n} = \sum_{m=0}^{n-1}[\hat{f}_{S_m^Y}(X_{m+1}) -
K_{S_m^Y}(\hat{f}_{S_m^Y})(X_m)]
\]
is such that $\{M_n,\mathcal{G}_n^X\}$ will be a martingale conditional
upon $\mathcal{G}_{\infty}^Y$. In addition, critical to our analysis,
will be that, $\mathbb{Q}^Y$-a.s., $\{S_n^Y(V)\}$ is finite. This
latter fact will enable us to control the various terms in
\eqref{eq:decomp} on events of the type $\{\sup_{k\geq0}S_k^Y(V)\leq
\bar{M}\}$ for $\bar{M}>0$. This is now elaborated.



\subsection{$\{M_m\}$ is $\mathbb{L}_p$-bounded}\label{63}

One can establish the following uniform in time $\mathbb{L}_p$-bounds
of the solution to Poisson's equation and the sequence $\{M_n\}$,
restricted to events $\{\sup_{k\geq0}S_k^Y(V)\leq\bar{M}\}$ for any
$\bar{M}>0$.

\begin{proposition} \label{regpoi1}
Assume (\textup{A1}). Let $r\in[0,1]$, $p\in[1,1/r]$
for $r \neq0$ and $p \geq1$ otherwise and $\bar{M}\in(0,\infty)$.
Then there exists $M<\infty$
such that for any $f\in\mathscr{L}_{V^r}$, $(x,y)\in E\times E$ and any
$m\in\mathbb{N}_0$,
\[
\mathbb{E}_{(x,y)}\bigl[|\hat{f}_{S_m^Y}(X_{m+1})|^p\mathbb{I}_{\{\sup
_{k\geq
0}S_k^Y(V)\leq\bar{M}\}}\bigr]^{1/p} \leq M V(x)^r.
\]
\end{proposition}

\begin{proposition} \label{martingale}
Assume (\textup{A1}). Let $r\in[0,1]$, $p\in[1,1/r]$
for $r \neq0$ and $p \geq1$ otherwise and $\bar{M}\in(0,\infty)$.
Then there exists $M<\infty$ such that for any $f\in\mathscr{L}_{V^r}$,
$(x,y)\in E\times E$ and any $m\in\mathbb{N}_0$,
\[
\mathbb{E}_{(x,y)}\bigl[ |M_{m}|^p\mathbb{I}_{\{\sup_{k\geq
0}S_k^Y(V)\leq\bar{M}\}}\bigr]^{1/p}  \leq m^{1/2\vee
1/p} M V(x)^r.
\]
\end{proposition}

This result will allow us to prove the $\mathbb{P}_{(x,y)}$-a.s.
convergence of $M_n$ to zero (cf. Appendix~B).

\subsection{Smoothness of the solution to Poisson's equation and $\omega(S_n^Y)$}\label{64}

As can be observed in \eqref{eq:decomp}, we have to control
the
fluctuations of the solution of the Poisson equation
$\{\hat{f}_{S_{m+1}^Y}(X_{m+1}) - \hat{f}_{S_{m}^Y}(X_{m+1})\}$. Also,
in \eqref{eq:prfdecomp}, the convergence of $\omega(S_m^Y)(f)$ to
$\omega(\eta)(f)$ $\mathbb{Q}_y$-a.s.~must be established. Both of
these issues are now dealt with.

\begin{proposition}\label{poisf1}
Assume \textup{(A1)} and \textup{(A2)}. Let
$r\in[0, 1)$, then for any $f\in\mathscr{L}_{V^{r}}$, $(x,y)\in
E\times
E$
\[
\lim_{m\rightarrow\infty}|\hat{f}_{S_{m+1}^Y}(X_{m+1}) -
\hat{f}_{S_m^Y}(X_{m+1})|=0\qquad \mathbb{P}_{(x,y)}\mbox{-a.s.}
\]
\end{proposition}

\begin{proposition}\label{prop:continuityomegaSnY}
Assume \textup{(A1)--(A3)}. Let $f\in\mathscr{L}_{V}$ and $(x,y)\in
E\times E$, then
\[
\lim_{m\rightarrow\infty} \omega(S_m^Y)(f)=\omega(\eta)(f)
\qquad\mathbb{Q}_{y}\mbox{-a.s.}
\]
\end{proposition}

\section{Examples}\label{7}

In this section we present some applications of our algorithms.
Specifically, it is demonstrated that the assumptions hold in some
very general scenarios. In addition, a numerical investigation of our
approach for a toy problem is given.

\subsection{Verifying the assumptions}\label{71}

It is now shown that it is possible to verify the assumptions in
Section \ref{sec:assumptions} in quite general scenarios. Let us
concentrate upon the case where, for $k\geq1$,
$(E,\mathcal{E})=(\mathbb{R}^k,\mathscr{B}(\mathbb{R}^k))$ and~$K$
(resp.,~$P$ -- recall the invariant measure is $\eta$) is a
symmetric random walk Metropolis kernel:
%
\begin{equation}\label{eq:rwmetropolis}
K(x,\mathrm{d}x') = \alpha_{\pi}(x,x')q_{\pi}(x-x')\,\mathrm{d}x' +
\delta_x(\mathrm{d}x')\biggl\{1-\int_{\mathbb{R}^k}
\alpha_{\pi}(x,x')q_{\pi}(x-x')\,\mathrm{d}x'\biggr\},
\end{equation}
where (resp.,~$P$)
\[
\alpha_{\pi}(x,x')  =  1\wedge\frac{\pi(x')}{\pi(x)}
\]
and $q_{\pi}$ (resp.~$q_{\eta}$) is a symmetric density (w.r.t.
Lebesgue measure).

\subsubsection{Assumptions}\label{711}

A set of general conditions is introduced, such that the assumptions in
Section \ref{sec:assumptions} will hold.

\begin{enumerate}[(M6)]
\item[(M1)]\textit{Density $\pi$.}
\begin{itemize}
\item$\pi$ admits a positive and continuous density w.r.t. Lebesgue
measure.
\end{itemize}
\item[(M2)]\textit{Definition of $\eta$.}
\begin{itemize}
\item$\eta(x)\propto\pi(x)^{\tilde{\alpha}}$, with
$\tilde{\alpha}\in(0,1)$.
\end{itemize}
\item[(M3)]\textit{Boundedness.}
\begin{itemize}
\item$\pi$ is upper bounded and bounded away from $0$ on compact sets.
\end{itemize}
\item[(M4)]\textit{Super-exponential densities.}
\begin{itemize}
\item$\pi$ is super exponential:
\[
\lim_{|x|\rightarrow+\infty}\frac{x}{|x|}\cdot\nabla\log(\pi(x)) =
-\infty.
\]
\end{itemize}
\item[(M5)]\textit{Regularity of contours.}
\begin{itemize}
\item The contours of $\pi$
are asymptotically regular:
\[
\limsup_{|x|\rightarrow+\infty}\frac{x}{|x|}\cdot\frac{\nabla\pi
(x)}{|\nabla\pi(x)|}<0.
\]
\end{itemize}
\item[(M6)]\textit{Lower bounds on $q_{\pi}$, $q_{\eta}$.}
\begin{itemize}
\item{Both $q_{\pi}$ and $q_{\eta}$ are such that there exists
$\tilde{\delta}_{q_{\pi}}>0$
(resp.,~$\tilde{\delta}_{q_{\eta}}>0$) and $\epsilon_{q_{\pi}}>0$
(resp.,~$\epsilon_{q_{\eta}}>0$) such that
\[
q_{\pi}(x) \geq\epsilon_{q_{\pi}}
\qquad\mbox{for }|x|<\tilde{\delta}_{q_{\pi}}
\]
(resp., $q_{\eta}(x) \geq\epsilon_{q_{\eta}}
\mbox{ for }|x|<\tilde{\delta}_{q_{\eta}}$). }
\end{itemize}
\end{enumerate}

\subsection{Result}\label{72}

\begin{proposition} \label{prop:verify}
Assume \textup{(M1)--(M6)}, then \textup{(A1)--(A3)} hold for any $r^*\in(0,1)$ with
\begin{eqnarray*}
W(x)
& = &
\biggl[\frac{|\pi|_{\infty}}{\pi(x)}\biggr]^{\tilde{\alpha}s_w},\qquad s_w\in(0,1),
\\
V(x)
& = &
\biggl[\frac{|\pi|_{\infty}}{\pi(x)}\biggr]^{s_v},\qquad s_v\in(0,r^*\tilde{\alpha} s_w).
\end{eqnarray*}
\end{proposition}

The proof is in Appendix \ref{sec:verifyproof}.

\subsubsection{Some comments}\label{721}

The conditions presented above are quite general. For example, they are
satisfied if $\pi$ is a mixture of normals. More generally, it may be
difficult to check the assumptions, but this is due to the underlying
nature of the geometric ergodicity assumptions; see \cite{jarner} for
more thorough investigations.

\subsection{Toy example}\label{sec:toyexample}

Our target distribution is
\[
\pi(x) = 0.4\psi(x;0,0.5) + 0.6\psi(x;17.5,1)
\]
with $\psi(x;\mu,\sigma^2)$ the normal density of mean $\mu$ and
variance $\sigma^2$.

Our algorithms are run with $K$ as a random walk Metropolis kernel with
normal random walk proposal density. The kernel is iterated 500 times
(i.e.,~$K=\widetilde{K}$ with $\widetilde{K}$ as a random walk
Metropolis kernel); this is to reduce the amount of interaction,
especially for large $\epsilon$. $\eta$ was taken to be:
\[
\eta(x) \propto\pi(x)^{0.75}.
\]
The algorithms were run for the same CPU time and the results can be
found, for 50 runs of the algorithm, in Table \ref{tab:nlres}. The
assumptions (M1)--(M6) are satisfied here.

\begin{table}[b]
\caption{Estimates from mixture comparison for
nonlinear MCMC. The estimates are for the expectation of~$X$; the true
value is 10.5. Each algorithm is run 50 times for 2 million iterations
after a 50~000 iteration burn-in (Section~\protect\ref{sec:auxselfinteration};
the simulations for Section~\protect\ref{sec:selection_mut} are adjusted for
the appropriate CPU time). The brackets are $\pm$2 standard
deviations across the repeats} \label{tab:nlres}
\begin{tabular*}{\textwidth}{@{\extracolsep{\fill}}llllll@{}}
\hline
Example & $\epsilon=0.05$ & $\epsilon=0.25$ & $\epsilon=0.5$ &$\epsilon=0.75$ & $\epsilon=0.95$\\
\hline
Section \ref{sec:selection_mut} & 10.32 ($\pm$0.08) & 10.74 ($\pm$0.12) & 10.89 ($\pm$0.19) & 10.37 ($\pm$0.18) & 10.99 ($\pm$0.20) \\
Section \ref{sec:auxselfinteration} & 10.57 ($\pm$0.04) &10.52 ($\pm$0.09) & 10.96 ($\pm$0.7) & 10.02 ($\pm$0.93) & 11.08 ($\pm$1.20)\\
\hline
\end{tabular*}
\end{table}

In Table \ref{tab:nlres}, the algorithms in Sections
\ref{sec:selection_mut} and \ref{sec:auxselfinteration} both perform
reasonably well for small values of $\epsilon$. As expected, from the
assumptions, as $\epsilon$ gets larger the accuracy falls. This is due
to the fact that the amount of auxiliary information that can enter
into the $\{X_n\}$ process is increased. For small $\epsilon$, the
example in Section \ref{sec:auxselfinteration} appears to work better
(more accurate estimation) due to the more sophisticated interaction
with the auxiliary chain. The drastic poor performance for the kernel
in Section \ref{sec:auxselfinteration}, for large $\epsilon$, is due to
the fact that no transition occurs after the swapping move.

To compare to the results of \cite{atchade2}, we ran a random walk
algorithm for 1 million iterations 50 times and a nonlinear algorithm
(Section~\ref{sec:auxselfinteration}). The nonlinear algorithm was run
with $\epsilon=0.01$ but the random walk kernel was not iterated. The
auxiliary chain was run with $\tilde{\alpha}=0.75$ (as in
(M2)). This was run for 110~000 iterations 50 times (which
is approximately the same CPU time as for the random walk Metropolis
algorithm). Both algorithms are such that all initial values are drawn
from a uniform on $[0,10.5]$. The estimated value for the first moment
is $6.93\pm16.96$ ($\pm $2 standard deviations, across the 50 runs)
and $10.41\pm2.03$ for the random walk and nonlinear methods,
respectively. The random walk algorithm is unable to jump between the
modes of the target, while the auxiliary chain is able to do so; hence
justifying our earlier intuition. This slightly contradicts the
`cautionary tale' in \cite{atchade2} as it illustrates that such
algorithms are potentially useful in cases where random walk algorithms
do not work well. We remark however, that one must be careful with
allowing too much auxiliary information to enter the chain
$\{X_n\}_{n\geq0}$; this can lead to poor results. This is consistent
with Proposition \ref{prop:K_mu_is_geometric}, which indicates that $d$
grows as $\epsilon$ goes to 1.

\section{Summary}\label{sec:summary}

We have investigated a new approach to stochastic simulation: Nonlinear
MCMC via auxiliary/self-interacting approximations. Convergence results for
several algorithms were established and the algorithm was demonstrated
on a toy example. As extensions to our ideas, the following may be
considered.

First, the conditions required for convergence may be relaxed. For
example, \cite{glynn} establishes weaker-than-geometric ergodicity
assumptions for the solution to the Poisson equation and functional
central limit theorem; also, \cite{fort} establishes drift conditions
for polynomial ergodicity. It would be of interest to see whether such
conditions would be sufficient for the convergence of our algorithms;
see \cite{roberts1} for proofs for parametric adaptive MCMC.

Second, it would be interesting to design more elaborate methods to
control the evolution of the empirical measure. In our current
algorithms, the empirical measure is only updated through the addition
of simulated points. It may enhance the algorithm to introduce some
mechanisms allowing the improvement of this quantity; for example, we
could introduce a death process with a rate associated with the
un-normalized target distribution.

\begin{appendix}\label{a1}

\section{Common properties of $K_{\mu}$}\label{app:commonpropkmu}

\begin{pf*}{Proof of Proposition \ref{prop:K_mu_is_geometric}}
The second and third statement of the proposition are a~direct
consequence of the first point from \cite{meyn1}, Theorem 2.3 (note the
$\phi$-irreducibility and aperiodicity follow immediately). The
minorization property is direct from the expression for $K_\mu$
and (A1)(ii) with
$\theta_{d}'=(1-\epsilon)\times
\theta_{d}$. Let us focus on the drift condition.

For any $x\in E$, $\mu\in\mathscr{P}_{\bar{M}}(E)$:
\[
K_{\mu}(V)(x) \leq(1-\epsilon)[\lambda V(x) + b \mathbb
{I}_{C_d}(x)] +
\epsilon[\mu(V) + V(x)\varphi(x) ],
\]
where $\varphi(x) = 1-\int_E \alpha(x,y)\mu(\mathrm{d}y)$. Then as
$\mu(V)<\bar{M}$, one has
\[
K_{\mu}(V)(x) \leq(1-\epsilon)[\lambda V(x) + b \mathbb
{I}_{C_d}(x)] +
\epsilon[\bar{M} + V(x)].
\]
Let $x\in C_d^c$, then
\[
K_{\mu}(V)(x) \leq\biggl[(1-\epsilon)\lambda+ \epsilon+
\frac{\epsilon\bar{M}}{d}\biggr]V(x)= \tilde{\lambda} V(x).
\]
For $x\in C_d$
\[
K_{\mu}(V)(x) \leq(1-\epsilon)[\lambda d + b] + \epsilon[\bar{M} + d]
\]
and hence one concludes that
\[
K_{\mu}(V)(x) \leq\tilde{\lambda} V(x) + \tilde{b}\mathbb{I}_{C_d}(x).
\]
\upqed\end{pf*}

\begin{pf*}{Proof of Proposition \ref{prop:invmeasnl3}}
This is a direct application of Proposition
\ref{prop:K_mu_is_geometric} and~Lem\-ma~%
\ref{lem:general_bound_variation_iterates_of_P_by_ variation_P}.
\end{pf*}

\begin{pf*}{Proof of Proposition \ref{prop:lipschitznl3}}
The proof is given for $r=1$ only. Let $|f|\leq V$:
\[
|[Q_{\mu}-Q_{\xi}](f)(x)| =
\biggl|\int_{E}[\mu-\xi](\mathrm{d}u)[\alpha(x,u)\{f(u)-f(x)\}]\biggr|.
\]
Now it is clear that, for any fixed $x\in E$:
\[
|\alpha(x,u)\{f(u)-f(x)\}| \leq[V(u) + V(x)],
\]
i.e.,
\[
|\alpha(x,u)\{f(u)-f(x)\}| \leq2V(u)V(x).
\]
Thus
\[
|[Q_{\mu}-Q_{\xi}](f)(x)| \leq 2V(x)\|\mu-\xi\|_V
\]
and then the result easily follows.
\end{pf*}

\section{Proof of the main result} \label{app:proofmainresult}

\begin{pf*}{Proof of Theorem \ref{llnnl3}}
Let $r \in[0,1)$ and $f \in\mathcal{L}_{V^r}$. Recall the strategy of
the proof outlined in Section \ref{sec:strategy}, which relies on the
decomposition:
%
\begin{equation}
S_n^X(f) -\pi(f)  =  S_n^X(f) - S^{\omega}_n(f) + S^{\omega}_n(f)
-\pi(f)
\end{equation}
with
\begin{eqnarray*}
&&(n+1)[S_n^X-S_n^{\omega}](f)
\\
&&\quad= M_{n+1} +
\sum_{m=0}^{n}[\hat{f}_{S_{m+1}^Y}(X_{m+1})-\hat
{f}_{S_m^Y}(X_{m+1})] +
\hat{f}_{S_0^Y}(X_0) - \hat{f}_{S_{n+1}^Y}(X_{n+1}) ,
\end{eqnarray*}
where $\{M_n\}$ is a martingale conditional upon
$\mathcal{G}_{\infty}^Y$. Proving the almost sure convergence of
$[S_n^X-S_n^{\omega}](f)$ relies on classical arguments. For any $n
\geq1$, $\delta>0$ and $\bar{M}\in(0,\infty)$,
\begin{eqnarray*}
&&
\mathbb{P}_{(x,y)}\Bigl( \sup_{k \geq n}|[S_k^X-S_k^{\omega}](f)|> \delta\Bigr)
\\
&&\quad\leq
\mathbb{P}_{(x,y)}\Bigl(\sup_{k\geq n}|M_{k+1}/(k+1)|>\delta/3,\sup_{k\geq0}S_k^Y(V)<\bar{M}\Bigr)
\\
&&\qquad{}+
\mathbb{P}_{(x,y)}\Biggl(\sup_{k \geq n}\Biggl|\sum_{m=0}^{k}[\hat{f}_{S_{m+1}^Y}(X_{m+1})-\hat{f}_{S_m^Y}(X_{m+1})]\Biggr|\Bigl/(k+1)> \delta/3,
\sup_{k\geq0}S_k^Y(V)<\bar{M}\Biggr)
\\
&&\qquad{}+
\mathbb{P}_{(x,y)}\Bigl(\sup_{k\geq n} [|\hat{f}_{S_0^Y}(X_0)| + |\hat
{f}_{S_{k+1}^Y}(X_{k+1})|]/(k+1)>\delta/3, \sup_{k\geq
0}S_k^Y(V)<\bar{M}\Bigr)
\\
&&\qquad{}+
\mathbb{Q}_y\Bigl(\sup_{k\geq0}S_k^Y(V)\geq\bar{M}\Bigr) .
\end{eqnarray*}
Let $\varepsilon>0$. By assumption there exists $\bar{M}>0$ such that
$\mathbb{Q}_y(\sup_{k\geq0}S_k^Y(V)\geq\bar{M}) \leq
\varepsilon/4$. Now we consider the remaining terms on the right-hand
side of the above equation from bottom to top; it is proved that there
exists $n_0>0$ such that for any $n \geq n_0$ each of these terms is
less than $\varepsilon/4$. Let $p \in(1,1/r)$. By Proposition
\ref{regpoi1}, one can apply Markov's inequality and a Borel--Cantelli
argument to show that the term on the third line vanishes as $n
\rightarrow\infty$. By Proposition \ref{poisf1} and a Ces\'aro
argument one concludes that the term on the second line goes to zero as
$n \rightarrow\infty$. The term dependent on $\{M_n\}$ is dealt with
by using an adaptation of a Birnbaum--Marshall inequality (see
\cite{andrieu1}) for $p \in(1,1/r)$.

%
%

Controlling the bias term requires a more novel approach. Note that
\[
|S_n^{\omega}(f)-\pi(f)| = \frac{1}{n+1}\Biggl|\sum_{i=0}^n[\omega
(S_i^Y) -
\omega(\eta)](f)\Biggr| ,
\]
as $\omega(\eta)=\pi$ in our setup. In Proposition
\ref{prop:continuityomegaSnY} it is proved that under our assumptions
$[\omega(S_i^Y) - \omega(\eta)](f) \rightarrow0$ $\mathbb{Q}_y$-a.s.
as $i \rightarrow\infty$. We conclude by invoking a Ces\`aro average
argument.
\end{pf*}

\begin{pf*}{Proof of Proposition \ref{regpoi1}}
Let $\bar{M}\in(0,\infty)$. The proof begins by conditioning upon the
filtration $\mathcal{G}^Y_m$ generated by the auxiliary process
$\{Y_n\}$; then, using the uniform in $\mu\in\mathscr{P}_{\bar{M}}(E)$,
geometric ergodicity is proved in Proposition
\ref{prop:K_mu_is_geometric}. As a result, there exists an $M<\infty$
such that
\begin{eqnarray*}
&&\mathbb{E}_{(x,y)}\bigl[|\hat{f}_{S_m^Y}(X_{m+1})|^p\mathbb{I}_{\{\sup_{k\geq 0}S_k^Y(V)\leq\bar{M}\}}\bigr]^{1/p}
\\
&&\quad\leq
M\mathbb{E}_{(x,y)}\bigl[|V(X_{m+1})^r|^p\mathbb{I}_{\{\sup_{k\geq0}S_k^Y(V)\leq\bar{M}\}}\bigr] ^{1/p}
\\
&&\quad\leq
M V^r(x),
\end{eqnarray*}
where we have used Jensen and the uniform drift condition on the set
$\{\sup_{k\geq0}S_k^Y(V)\leq\bar{M}\}$ proved in Proposition
\ref{prop:K_mu_is_geometric}
\end{pf*}

\begin{pf*}{Proof of Proposition \ref{martingale}}
We follow a similar argument to that of \cite{andrieu1}, Proposition~6.
Throughout, denote by $B_p$ a generic constant dependent upon $p$ only.
Also recall $pr\leq1$. The proof begins by applying the
B\"urkholder--Davis inequality (see, e.g.,~\cite{shiryaev},
pages~499--500), which yields for $p \geq1$
\begin{eqnarray*}
&&
\mathbb{E}_{(x,y)}\bigl[|M_{n}|^p\mathbb{I}_{\{\sup_{k\geq0}S_k^Y(V)\leq\bar{M}\}}\bigr]^{1/p}
\\
&&\quad\leq
B_p\mathbb{E}_{y}\Biggl[\mathbb{E}_{(x,y)}\Biggl[\Biggl(\sum_{m=0}^{n-1}[\hat{f}_{S_m^Y}(X_{m+1})-K_{S_m^Y}(\hat{f}_{S_m^Y})(X_m)]^2\Biggr)^{p/2}\Bigl|\mathcal{G}^Y_{\infty}\Biggr]
\mathbb{I}_{\{\sup_{k\geq0}S_k^Y(V)\leq\bar{M}\}}\Biggr]^{1/p}.
\end{eqnarray*}
In the case $p>2$, by similar manipulations to those featured in
\cite{andrieu1}
\[
\mathbb{E}_{(x,y)}\bigl[ |M_{n}|^p\mathbb{I}_{\{\sup_{k\geq
0}S_k^Y(V)\leq\bar{M}\}}\bigr]^{1/p} \leq n^{1/2} B_p M V(x)^{r}.
\]
In the case $p \leq2$, one may apply the $C_p$-inequality to yield
\begin{eqnarray*}
&&
\mathbb{E}_{(x,y)}\bigl[ |M_{n}|^p\mathbb{I}_{\{\sup_{k\geq
0}S_k^Y(V)\leq\bar{M}\}}\bigr]^{1/p}
\\
&&\quad\leq
\Biggl[\sum_{m=0}^{n-1}\mathbb{E}_{(x,y)}\bigl[ |\hat{f}_{S_m^Y}(X_{m+1})
- K_{S_m^Y}(\hat{f}_{S_m^Y})(X_m)|^p \mathbb{I}_{\{\sup_{k\geq
0}S_k^Y(V)\leq\bar{M}} \bigr]\Biggr]^{1/p}.
\end{eqnarray*}
Application of Minkowski, conditional Jensen and Proposition
\ref{regpoi1} yields
\[
\mathbb{E}_{(x,y)}\bigl[ |M_{n}|^p\mathbb{I}_{\{\sup_{k\geq
0}S_k^Y(V)\leq\bar{M}\}}\bigr]^{1/p} \leq n^{1/p}M V^r(x)
\]
from which we can conclude.
\end{pf*}

\begin{pf*}{Proof of Proposition \ref{poisf1}}
Our proof is based upon the decomposition of Proposition \ref{iter22}
(in Appendix \ref{app:standardMC}) and then using the Lipschitz
continuity properties proved in Propositions \ref{prop:invmeasnl3} and
\ref{prop:lipschitznl3}. Let $\bar{M}\in(0,\infty)$ be given; suppose
that we are on the set $\{\sup_{k\geq0}S_k^Y(V)\leq\bar{M}\}$. Then
%
\begin{eqnarray}\label{eq:poisdecomp}
&&|\hat{f}_{S_{m+1}^Y}(X_{m+1}) - \hat{f}_{S_m^Y}(X_{m+1})|\quad\nonumber
\\
&&\quad=
\Biggl|\sum_{n\in\mathbb{N}}\sum_{i=0}^{n-1}[K_{S_{m+1}^Y}^i-\omega
(S_{m+1}^Y)](K_{S_{m+1}^Y}
- K_{S_{m}^Y})[K_{S_m^Y}^{n-i-1}-\omega
(S_m^Y)(f)(X_{m+1})]\quad
\\
&&\qquad{}-
\sum_{n\in\mathbb{N}}[\omega(S_{m+1}^Y)-\omega
(S_m^Y)]\bigl(K_{S_m^Y}^n-\omega(S_m^Y)\bigr)(f)\Biggr|.\quad\nonumber
\end{eqnarray}

Now, consider the first term. Since, for any $m\geq0$, the kernel
$K_{S_m^Y}$ satisfies:
\[
\|[K_{S_m^Y}^n - \omega(S_m^Y)](f)\|_{V^r} \leq M\rho^n V(X_{m+1})^r
\]
for some finite $M$ and $\rho\in(0,1)$ independent of $S_m^Y \in
\mathscr{P}_{\bar{M}}(E)$, it follows that:
\begin{eqnarray*}
&&|[K_{S_{m+1}^Y}^i-\omega(S_{m+1}^Y)](K_{S_{m+1}^Y} -
K_{S_{m}^Y})[K_{S_m^Y}^{n-i-1}-\omega(S_m^Y)(f)(X_{m+1})]| \\
&&\quad\leq
M\rho^{i}V(X_{m+1})^r |(K_{S_{m+1}^Y} -
K_{S_{m}^Y})[K_{S_m^Y}^{n-i-1}-\omega(S_m^Y)(f)]|_{V^r}.
\end{eqnarray*}
Then, adopting the continuity result for $K_{S^Y_m}$:
\[
\|\!|K_{\mu}-K_{\xi}|\!\|_{V^r} \leq2\|\mu-\xi\|_{V^r}
\]
for \textit{any} $\mu,\xi\in\mathscr{P}_{\infty}(E)$, it follows that:
\[
|(K_{S_{m+1}^Y} -
K_{S_{m}^Y})[K_{S_m^Y}^{n-i-1}-\omega(S_m^Y)(f)]|_{V^r} \leq
M\rho^{n-i-1}\|S_{m+1}^Y-S_m^Y\|_{V^r}.
\]
Since $\|S_{m+1}^Y-S_m^Y\|_{V^r}\leq[V(Y_{m+1})^r + S_m^Y(V^r)]/(m+2)$
\begin{eqnarray*}
&&\sum_{n,i}|[K_{S_{m+1}^Y}^i-\omega(S_{m+1}^Y)](K_{S_{m+1}^Y}
- K_{S_{m}^Y})[K_{S_m^Y}^{n-i-1}-\omega(S_m^Y)(f)(X_{m+1})]|
\\
&&\quad\leq
\frac{M}{(1-\rho)^2}\frac{V(X_{m+1})^r}{m+2}[V(Y_{m+1})^r +
S_m^Y(V^r)].
\end{eqnarray*}

Turning to the second sum on the right-hand side of
(\ref{eq:poisdecomp}), using the continuity result
\[
\|\omega(\mu)-\omega(\xi)\|_{V^r} \leq M\Vert\!\vert K_{\mu}-K_{\xi}\vert\!\Vert_{V^r}
\]
(for $M<\infty$ not depending on $\mu,\xi\in\mathscr{P}_{\bar{M}}(E)$
by Proposition \ref{prop:lipschitznl3}) and the continuity of the
kernel $K_{\mu}$ (Lemma
\ref{lem:general_bound_variation_iterates_of_P_by_ variation_P} )
yields,
\[
\bigl|[\omega(S_{m+1}^Y)-\omega(S_m^Y)]\bigl(K_{S_m^Y}^n-\omega(S_m^Y)\bigr)(f)\bigr|
\leq
M\rho^n \frac{[V(Y_{m+1})^r + S_m^{Y}(V^r)]}{m+2},
\]
from which we obtain a similar bound for the second sum on the
right-hand side of (\ref{eq:poisdecomp}).

We now establish an $\mathbb{L}_p$-bound, for $p>1$ of this upper bound
on $\{\sup_{k\geq0}S_k^Y(V)\leq\bar{M}\}$, which will allow us to use
a Borel--Cantelli argument to complete the proof. Note that it is
naturally sufficient to consider $\frac{V(X_{m+1})^r}{m+2}V(Y_{m+1})^r$
on $\{\sup_{k\geq0}S_k^Y(V)\leq\bar{M}\}$, and we focus on
%
\begin{eqnarray}
&&
\mathbb{E}_{(x,y)}\bigl[V(X_{m+1})^{rp}V(Y_{m+1})^{rp}\mathbb{I}_{\{\sup_{k\geq0}S_k^Y(V)\leq\bar{M}\}}\bigr]^{1/p}\nonumber
\\
&&\quad=
\mathbb{E}_{(x,y)}\bigl[\mathbb{E}_{(x,y)}[V(X_{m+1})^{rp}|\mathcal{G}^Y_\infty] V(Y_{m+1})^{rp} \mathbb{I}_{\{\sup_{k\geq0}S_k^Y(V)\leq\bar{M}\}}\bigr]^{1/p}
\\
&&\quad\leq
M V(x)^{r} \mathbb{E}_{y}\bigl[V(Y_{m+1})^{rp} \mathbb{I}_{\{\sup_{k\geq0}S_k^Y(V)\leq\bar{M}\}}\bigr]^{1/p}\leq M V(x)^r W^{r r^*}(y),\nonumber
\end{eqnarray}
where we have used that, conditional upon $\mathcal{G}^Y_\infty$ and on
the event $\{\sup_{k\geq0}S_k^Y(V)\leq\bar{M}\}$, the following bound
holds
$\mathbb{E}_{(x,y)}[V(X_{m+1})^{pr}|\mathcal{G}^Y_{\infty
}]^{1/p}\leq
M_0 V(x)^r$ for some deterministic constant $M_0$ depending only on
$\bar{M}$ and the parameters of the drift condition in (A1). Similarly,
$M \geq M_0$ only depends on $\bar{M}$ and the parameters of the drift
conditions in (A1)--(A2). With $p>1$ we conclude that
\[
\sum_{m=0}^{\infty}\mathbb{P}_{(x,y)}\biggl(\biggl\{\frac
{1}{m+1}V(X_{m+1})^{r}[V(Y_{m+1})^r
+ S_m^Y(V^r)] > \varepsilon, \sup_{k\geq0}S_k^Y(V)\leq
\bar{M}\biggr\}\biggr) < \infty.
\]
The result then follows by using for any $\delta>0$ the bound,
\begin{eqnarray*}
&&
\mathbb{P}_{(x,y)}\Bigl(\sup_{k\geq m} |\hat{f}_{S_{k+1}^Y}(X_{k+1}) -
\hat{f}_{S_k^Y}(X_{k+1})|>\delta\Bigr)
\\
&&\quad\leq
\mathbb{Q}_y\Bigl(\sup_{k\geq0}S_k^Y(V)\geq
\bar{M}\Bigr)
+
 \mathbb{P}_{(x,y)}\Bigl(\sup_{k\geq m} |\hat{f}_{S_{k+1}^Y}(X_{k+1})
- \hat{f}_{S_k^Y}(X_{k+1})|>\delta, \sup_{k\geq0}S_k^Y(V)\leq
\bar{M}\Bigr)
\end{eqnarray*}
and using the fact that for any $\varepsilon>0$ one can find an
$\bar{M}$ large enough to ensure that the first term on the right-hand
side is less than $\varepsilon/2$ and then $m_0$ such that for any $m
\geq m_0$ the second term on the right-hand side is also upper bounded
by $\varepsilon/2$.
\end{pf*}

\begin{pf*}{Proof of Proposition \ref{prop:continuityomegaSnY}}
Note first that for any $i,j\in\mathbb{N}$, $f \in\mathcal{L}_{V}$ and
$x \in E$ such that $V(x)<\infty$ we have the following bound,
\begin{eqnarray*}
&&
|[\omega(S_i^Y) - \omega(\eta)](f)|
\\
&&\quad\leq
|\omega(S_i^Y)(f) -K_{S_i^Y}^j(f)(x)|+ |K_{S_i^Y}^j(f)(x)- K_{\eta}^j(f)(x)| +|K_{\eta}^j(f)(x)-\omega(\eta)(f)| .
\end{eqnarray*}
Let $\varepsilon,\delta>0$ and $\bar{M}>\eta(V)$ be such that
$\mathbb{Q}_y(\sup_{k \geq0}S^Y_k(V) \geq\bar{M})<\varepsilon/4$. On
the event $\{\sup_{k \geq0}S^Y_k(V) < \bar{M}\}$ we have by
Proposition \ref{prop:invmeasnl3} the existence of $M<+\infty$ and
$\rho\in[0,1)$ (independent of $i$) such that the first and last
terms on the right-hand side are bounded by $M\rho^j$. We can therefore
fix $m$ such that
\[
\mathbb{Q}_{y}\Bigl(\sup_{j \geq m, i \geq0}|\omega(S_i^Y)(f) -
K_{S_i^Y}^j(f)(x)| + |K_{\eta}^j(f)(x)-\omega(\eta)(f)| > \delta/2,
\sup_{k \geq0}S^Y_k(V) < \bar{M}\Bigr) \leq\varepsilon/2 .
\]
Now from Lemma \ref{lemma:kiterates} one may conclude that there exists
$m_0>0$ such that for any $m \geq m_0$
\[
\mathbb{Q}_{y}\Bigl(\sup_{i \geq m}|K_{S_i^Y}^j(f)(x)-
K_{\eta}^j(f)(x)| > \delta/2, \sup_{k \geq0}S^Y_k(V) < \bar{M}\Bigr)
\leq\varepsilon/4 .
\]
The proof is completed by noting that the results above imply that for
$m \geq m_0$,
%
\begin{equation}
\mathbb{Q}_{y}\Bigl( \sup_{i \geq m}|[\omega(S_i^Y) -
\omega(\eta)](f)|>\delta\Bigr) \leq\varepsilon.
\vspace*{5pt}\end{equation}
\upqed\end{pf*}

\section{Standard technical results on Markov chains} \label{app:standardMC}

\begin{lemma}\label{lem:general_bound_variation_iterates_of_P_by_ variation_P}
Let $(E,\mathcal{E})$ be a measurable space, $\bar{b}<\infty$,
$\bar{\lambda} \in(0,1)$ and $\bar{C}\in\mathcal{E}$. Then for any
Markov transition probabilities $P_1,P_2 \dvtx E \rightarrow
\mathscr{P}(E)$ satisfying for $(x,A)\in E \times\mathcal{E}$ and
$i=1,2$,
%
\begin{eqnarray}
P_iV(x) &\leq&\bar{\lambda}V(x)+\mathbb{I}_{\bar{C}}(x)\bar{b} ,
\\
P_i(x,A) &\geq&\mathbb{I}_{\bar{C}}(x)\bar{\epsilon} \bar{\nu}(A).
\end{eqnarray}
There exist $\bar{M}(\cdot)<\infty$, $\bar{\rho} \in[0,1)$, invariant
probability measures $\pi_1,\pi_2 \in\mathscr{P}(E)$ (corresponding to
$P_1$ and $P_2$, respectively), such that for any $n \geq1$,
$r\in[0,1]$ and any $|f|\leq V^r$
\[
|[P_1^n - \pi_1](f)|_{V^r} \vee|[P_2^n - \pi_2](f)|_{V^r} \leq
\bar{M}(r) \bar{\rho}^n
\]
for any $n \geq1$,
\[
\Vert\!\vert P^n_1-P^n_2\vert\!\Vert_{V^r} \leq \bar{M}(r)\Vert\!\vert P_1 - P_2\vert\!\Vert_{V^r}
\]
and
\[
\Vert \pi_1-\pi_2\Vert_{V^r}  \leq \bar{M}(r)\Vert\!\vert P_1 - P_2\vert\!\Vert_{V^r} .
\]
\end{lemma}

\begin{pf}
Let $r\in[0,1]$ and $f \in\mathscr{L}_{V^r}$.
We have the following decomposition:
%
\[
|[P_{1}^n-P_{2}^n](f)|  = \Biggl|\sum_{i=0}^{n-1} P_{1}^i\bigl([P_{1}-
P_{2}]\{ [P_{2}^{n-i-1}-\pi_2](f)\}\bigr)\Biggr|.
\]

For any $|f|\leq V^r$, in a similar manner to Proposition 3 of
\cite{andrieu1}:
\begin{eqnarray*}
|[P_{1}^n-P_{2}^n](f)|
& \leq&
\bar{M}(r)\sum_{i=0}^{n-1}
\bar{\rho}^{n-i-1}
P_{1}^i(\Vert P_{1}- P_{2}\Vert_{V^r})
\\
& = &
\bar{M}\sum_{i=0}^{n-1} \bar{\rho}^{n-i-1}
P_{1}^i\biggl(\frac{\Vert P_{1}- P_{2}\Vert_{V^r}}{V^r} V^r\biggr)
\\
& \leq&
\bar{M}(r)\Vert\!\vert P_{1}- P_{2}\vert\!\Vert_{V^r}\sum_{i=0}^{n-1}
\bar{\rho}^{n-i-1} P_{1}^i( V^r).
\end{eqnarray*}
From the drift condition (A2) and conditional Jensen one can bound
$P^i_1 V^r$ by $[\bar{\lambda} + \bar{b}/(1-\bar{\lambda})]^{r}V(x)^r$
for $r \in[0,1]$ and hence conclude that:
\[
|[P_1^n-P_{2}^n](f)|  \leq \bar{M}(r)\Vert\!\vert P_{1}- P_{2}\vert\!\Vert_V^r.
\]
Since the right-hand side is independent of $n$, the inequality holds
in the limit and hence, by $V$-uniform ergodicity, the result.
\end{pf}

\begin{proposition}\label{iter22}
Assume
(\textup{A1}). 
Then, for $r\in[0,1]$, $\xi,\mu\in\mathscr{P}_{\infty}(E)$, $f
\in
\mathscr{L}_{V^r}$ we have the following decomposition for the
differences in the solution to the Poisson equation:
\begin{eqnarray*}
\hat{f}_\xi(x) - \hat{f}_\mu(x)
& = &
\sum_{n\in\mathbb{N}}\Biggl\{\sum_{i=0}^{n-1}\bigl([K^i_{\xi}-\omega(\xi)](K_{\xi}-K_{\mu})\{[K^{n-i-1}_{\mu}-\omega(\mu)](f)\}(x)\bigr)
\\
&&\hphantom{\sum_{n\in\mathbb{N}}\Biggl\{}{}-
[\omega(\xi)-\omega(\mu)]\bigl([K^n_{\mu}-\omega(\mu)](f)\bigr) \Biggr\}.
\end{eqnarray*}
\end{proposition}

\begin{pf}
Adopting the resolvent solution to the Poisson equation (which exists
under our assumptions), we have
\begin{eqnarray*}
\hat{f}_\xi(x) - \hat{f}_\mu(x)
& = &
\sum_{n\in\mathbb{N}_0}\bigl[\bigl([K_{\xi}^n-\omega(\xi)](f)(x)\bigr) - \bigl([K_{\mu}^n-\omega(\mu)](f)(x)\bigr)\bigr]
\\
& = &
\sum_{n\in\mathbb{N}}\Biggl[\sum_{i=0}^{n-1} K_{\xi}^i\bigl([K_{\xi}-K_{\mu}]\{[K_{\mu}^{n-i-1}-\omega(\mu)](f)\}\bigr)(x)+\omega(\mu)(f)-\omega(\xi)(f)\Biggr]
\\
& = &
\sum_{n\in\mathbb{N}}\Biggl\{\sum_{i=0}^{n-1}\bigl([K^i_{\xi}-\omega(\xi)](K_{\xi}-K_{\mu})\{[K^{n-i-1}_{\mu}-\omega(\mu)](f)\}(x)\bigr)
\\
&&\hphantom{\sum_{n\in\mathbb{N}}\Biggl\{}{}-
[\omega(\xi)-\omega(\mu)]\bigl([K^n_{\mu}-\omega(\mu)](f)\bigr)\Biggr\}
\end{eqnarray*}
since
\[
-\sum_{i=0}^{n-1}\omega(\xi)[K_{\xi}-K_{\mu}](K_{\mu}^{n-i-1}(f))
 =
-\omega(\xi)\bigl(f-K_{\mu}^n(f)\bigr).
\]
\end{pf}

\section{Convergence of the iterates}

The main result of this section is Lemma \ref{lemma:kiterates}, where
it is established that for any $q\geq1$, $f\in\mathscr{L}_{V}$
%
\begin{equation}\label{eq:appendixcmainresult}
\lim_{n\rightarrow\infty}|[K_{S_n^Y}^q - K_{\eta}^q](f)(x)| = 0,
\qquad\mathbb{Q}_y\mbox{-a.s.},
\end{equation}
with $K_{\mu}$ as in \eqref{eq:generic_K_mu}. The proof consists of
showing that $K_{\mu}^q(f)(x)$ can be rewritten as $\mu^{\otimes
q}(\overline{g})$ for some function $\overline{g}\dvtx E^q \rightarrow
\mathbb{R}$ to be given below. We will then use results from Appendix
\ref{appendix:Vstatistic}, associated with $V$-statistics for an
appropriate class of functions, to complete our argument.

Introduce the following family of Markov transition probabilities on
$(E \times E,\mathcal{E}\otimes\mathcal{E})$, indexed by $z_1 \in E$,\vadjust{\goodbreak}
\begin{eqnarray*}
&&
T_{z_1}((w_0, w_0');\mathrm{d}(w_1,w_1'))
\\
&&\quad:=
(1-\epsilon)
K(w_0,\mathrm{d}w_1) \delta_{w_0}(\mathrm{d}w_1')
\\
&&\qquad{}+
\epsilon\bigl[\alpha(w_0,z_1)\delta_{(z_1,w_0)}(\mathrm{d}w_1,\mathrm{d}w_1')+\bigl(1-\alpha(w_0,z_1)\bigr) \delta_{(w_0,w_0')}(\mathrm{d}w_1,\mathrm{d}w_1')\bigr].
\end{eqnarray*}
For any $w_0,w_0'\in E$ and $z:=(z_1,\ldots,z_q) \in E^q,$ we define
the iterates of this family of kernels as follows: for $k=2,\ldots,q$
and any $f\in\mathscr{L}_V$,
%
\begin{equation}\label{eq:titerates_appendixc}
T_{z_1,\ldots,z_k}^k(f\otimes1)(w_0, w_0') :=
T_{z_1,\ldots,z_{k-1}}^{k-1} \bigl(T_{z_k}(f\otimes1)(\cdot)\bigr)(w_0, w_0'),
\end{equation}
where for any $x,x' \in E$, $(f\otimes1)(x,x'):=f(x)$. Let
$z:=(z_1,\dots,z_q)\in E^q$. Following an argument identical to that
developed in the proof of Lemma \ref{lemma:kiterates} it is possible to
show that for any $k = 1, \ldots, q$ $T_{z_1,\ldots,z_k}^k(f\otimes
1)(w_0, w_0')$ belongs to $\mathscr{L}_{\mathcal{V}_{z_1,\ldots,z_k}}$
where for $w,w'\in E$,
\[
\mathcal{V}_{z_1,\ldots,z_k}(w,w'):= V(w)+V(w')+\sum_{i=1}^k V(z_i) .
\]

\begin{qiterates}\label{prop:qiterates_new}
Assume (\textup{A1}). For any $q\geq1$, $(z_1,\ldots,z_q) \in E^q$,
$\mu\in\mathscr{P}_\infty(E)$, $f\in\mathscr{L}_V$, $x,x'\in E$ we
have that
\[
K_{\mu}^q(f)(x) = \int_{E^{q}} T_{z_1,\ldots,z_q}^q(f\otimes1)(x, x')
\mu^{\otimes q} ( \mathrm{d}(z_1,\dots,z_q) ) .
\]
\end{qiterates}

\begin{pf}
The result is proved by induction. One immediately checks that for any
$z_1 \in E$, $f\in\mathscr{L}_V$, $w_0,w_0'\in E$,
\[
T_{z_1}(f \otimes1)(w_0, w_0') = (1-\epsilon) K(f)(w_0) + \epsilon
\bigl[\alpha(w_0,z_1) f(z_1) + \bigl(1-\alpha(w_0,z_1)\bigr) f(w_0)\bigr],
\]
and hence
\[
\mu\bigl(T_{z_1}(f \otimes1)(w_0, w_0')\bigr) = \int_E T_{z_1}(f
\otimes1)(w_0, w_0') \mu(\mathrm{d}z_1) = K_\mu(f)(w_0).
\]
Now assume the property is true for $k-1 \geq1$. Then
\begin{eqnarray*}
\mu^{\otimes k} \bigl( T_{z_1,\ldots,z_k}^k(f\otimes1)(w_0, w_0') \bigr)
&=&
\mu^{\otimes(k-1)}\bigl(T_{z_1,\ldots,z_{k-1}}^{k-1}\bigl\{\mu
\bigl(T_{z_k}(f\otimes1)(\cdot)\bigr)\bigr\}(w_0, w_0')\bigr)
\\
& = &
\mu^{\otimes(k-1)}\bigl(T_{z_1,\ldots,z_{k-1}}^{k-1}
\bigl(K_\mu(f)\otimes1\bigr)(w_0, w_0')\bigr) ,
\end{eqnarray*}
as required.
\end{pf}

Now, to establish \eqref{eq:appendixcmainresult} we need to show that
$T_{z_1,\ldots,z_q}^q(f)(w_0, w_0')$ lies within the class of functions
for which Lemma \ref{lemma:vstatistics} applies; this is proved below.

\begin{kiterates}\label{lemma:kiterates}
Assume \textup{(A1)--(A3)}. Let $q\geq1$ be fixed and $f\in\mathscr{L}_V$. Then
for any $x\in E$
\[
\lim_{n\rightarrow\infty}|[K_{S_n^Y}^q-K_{\eta}^q](f)(x)| =
0\qquad\mathbb{Q}_y\mbox{-a.s.}
\]
\end{kiterates}

\begin{pf}
Our objective is to use the representation established in Proposition
\ref{prop:qiterates_new} along with the result in Lemma
\ref{lemma:vstatistics}. To that end we show that for any
$f\in\mathscr{L}_{V}$, then $T_{z_1,\ldots,z_q}^q(f\otimes1)(w_0,
w_0')\in\mathscr{L}_{\mathcal{V}_z^{(q)}}$, $z^{(q)}=(z_1,\dots,z_q)$,
where $T_{z_1,\ldots,z_q}^q(f\otimes1)(w_0, w_0')$ is as in
\eqref{eq:titerates_appendixc}.
The result can be proved by induction. Now, for any $k=1,\ldots,q$,
$w_{k-1}, w_{k-1}' \in E$ and $z=(z_1,\ldots,z_q) \in E^q$
\begin{eqnarray*}
T_{z_k}\bigl(\mathcal{V}_{z^{(q)}}\bigr)(w_{k-1}, w_{k-1}')
& := &
(1-\epsilon)\Biggl[K(V)(w_{k-1}) + V(w_{k-1})+ \sum_{i=1}^q V(z_i)\Biggr]
\\
&&{}+
\epsilon\Biggl\{ \alpha(w_{k-1},z_k)\Biggl[V(w_{k-1}) + V(z_k) + \sum_{i=1}^q V(z_i)\Biggr]
\\
&&\hphantom{+\epsilon\Biggl\{}{}+
\bigl(1-\alpha(w_0,z_1)\bigr) \Biggl[ V(w_{k-1}) + V(w_{k-1}') + \sum_{i=1}^q
V(z_i) \Biggr]\Biggr\}.
\end{eqnarray*}
Since there exists $M<\infty$ such that for any $x \in E$, $K(V)(x)
\leq M V(x)$ we conclude that there exists $C_1>0$ such that for any
$k=1,\ldots,q$, $w_{k-1}, w_{k-1}' \in E$ and $z^{(q)}\in E^q$
%
\begin{equation}
T_{z_k}(\mathcal{V}_{z^{(q)}})(w_{k-1}, w_{k-1}') \leq C_1
\mathcal{V}_{z^{(q)}}(w_{k-1},w_{k-1}').
\end{equation}
This implies that for any $g \in\mathscr{L}_{\mathcal{V}_{z^{(q)}}}$
then $T_{z_k}(g)(w_{k-1}, w_{k-1}') \in
\mathscr{L}_{\mathcal{V}_{z^{(q)}}}$. Now we can proceed with the
induction. Assume that for some $k-1\geq1$, if $g \in
\mathscr{L}_{\mathcal{V}_{z^{(q)}}}$, then
$T^{k-1}_{z_{1},\ldots,z_{k-1}}(g)(w, w') \in
\mathscr{L}_{\mathcal{V}_{z^{(q)}}}$. Then by definition
\[
T_{z_1,\ldots,z_k}^k(f\otimes1)(w_0, w_0') =
T_{z_1,\ldots,z_{k-1}}^{k-1}\{T_{z_k}(f\otimes1)(\cdot)\}(w_0, w_0'),
\]
and the induction follows. Now, for any fixed $w_0,w_0'$ one has that
$T_{z_1,\ldots,z_q}^q(f\otimes1)(w_0, w_0')\in\mathscr{L}_{W^{(q)}}$
and the result follows from Lemma \ref{lemma:vstatistics}.
\end{pf}

\section{Results on $U$ and $V$-statistics for Markov chains} \label{appendix:Vstatistic}

Let $(E,\mathcal{E})$ be a Polish space and
$\eta\in\{\mu\in\mathscr{P}(E)\dvtx \mu(W)<\infty\}$. Denote $\Omega=
E^{\mathbb{N}}$ and $\mathcal{F} = \mathcal{E}^{\otimes\mathbb{N}}$
and consider a time-homogeneous Markov chain $\{X_n\}_{n\geq0}$ with
transition kernel $P$ such that $\eta P = \eta$ with $X_0=x$. Denote by
$\mathbb{P}_x$ the corresponding probability distribution. Note that
$\{X_n\}$ should not be confused with the process introduced in Section
\ref{sec:algo_descrip}.

For any sequence $\{Z_n\}$, $Z_n\in E$, any $q \in\mathbb{N}$ and
$f\dvtx E^q \rightarrow\mathbb{R}$, denote for any $n \geq1$ the
associated $V$-statistic
%
\begin{equation}
S_{n,Z}^{\otimes q}(f) = \frac{1}{(n+1)^q}
\sum_{\vartheta\in(q,n+1)}f\bigl(Z_{\vartheta(1)},\dots,Z_{\vartheta(q)}\bigr),
\end{equation}
where $(q,n+1)$ is the set of all mappings of $\{0,\dots,q-1\}$ into
$\{0,\dots,n\}$.

The main result of this section is Lemma \ref{lemma:vstatistics}, where
it is shown that under additional assumptions on $P$ and $f$, that
\[
\lim_{n\rightarrow\infty} S_{n,X}^{\otimes q}(f) = \eta^{\otimes q}(f),
\]
$\mathbb{P}_x$-a.s. The proof relies on a coupling argument with
another Markov chain $\{Y_n\}_{n\geq0}$ defined on
$(\Omega,\mathcal{F})$ with the same transition $P$, but initialized at
stationarity, that is,~\mbox{$Y_0 \sim\eta$}. $\mathbb{P}_\eta$ denotes the
corresponding probability distribution.

The conditions on $\{X_n\}_{n\geq0}$ and $\{Y_n\}_{n\geq0}$ referred
to above are given in (A2), and will, in
particular, imply geometric ergodicity. The class of functions to which
our results apply is defined as follows. Let $(W^r)^{(q)}(x^{(q)}):=
\sum_{i=1}^q W(x_i)^r$ for any $r\in(0,1)$, $x^{(q)}:=(x_1,\ldots,x_q)
\in E^q$; we will consider below the following class of functions
\[
\mathscr{L}_{ (W^r)^{(q)} } := \Bigl\{ f\in m E^q\dvtx \sup_{x^{(q)} \in
E^q} \bigl|f\bigl(x^{(q)}\bigr)\bigr|/(W^r)^{(q)}\bigl(x^{(q)}\bigr) < \infty\Bigr\}.
\]

For any sequence $\{Z_n\}$, $Z_n\in E$, any $q \in\mathbb{N}$ and
$f\dvtx E^q \rightarrow\mathbb{R}$ denote for any $n \geq1$ the associated
$U$-statistic
%
\begin{equation}\label{eq:ustat_new}
S_{n,Z}^{\odot q}(f) =
\frac{1}{(n+1)_q}\sum_{\vartheta\in\langle q,n+1\rangle}f\bigl(Z_{\vartheta
(1)},\dots,Z_{\vartheta(q)}\bigr),
\end{equation}
where $\langle q,n+1\rangle$ is the set of one-to-one mappings from
$\{0,\dots,q-1\}$ into $\{0,\dots,n\}$ and $n_q := n!/(n-q)!$. A
preliminary result on $U$-statistics is first established, based on the
aforementioned coupling.

\begin{coupu}\label{appendixprop:markovcoupl}
Assume \textup{(A2)} and \textup{(A3)}. Let $\{X_n\}_{n\geq0}$ and $\{Y_n\}_{n\geq0}$ be as
defined above. Then for any $q \in\mathbb{N}$, $r\in[0,1)$,
$f\in\mathcal{L}_{(W^r)^{(q)}}$ and $x\in E,$ there exists a coupling
$\{\check{X}_n,\check{Y}_n\}_{n\geq0}$ on some probability space
$(\Omega\times\Omega,\mathcal{F}\otimes\mathcal{F},\widetilde
{\mathbb{P}})$,
such that
\[
\lim_{n\rightarrow\infty}|S_{n,\check{X}}^{\odot q}(f) -
S_{n,\check{Y}}^{\odot q}(f)| = 0
\qquad\widetilde{\mathbb{P}}\mbox{-a.s.}
\]
\end{coupu}

\begin{pf}
Let $\mathbb{P}_{x}^{(n)}$ (resp.,~$\mathbb{P}_{\eta}^{(n)}$)
denote the law of $(X_n,X_{n+1},\dots)$
(resp.,~$(Y_n,Y_{n+1},\dots)$). Then, convergence in total
variation of the processes is sufficient to imply that:
\[
\lim_{n\rightarrow\infty}\bigl\|\mathbb{P}_x^{(n)}-\mathbb{P}_{\eta
}^{(n)}\bigr\|_{\mathrm{TV}}
 =  0.
\]
By Theorem 2.1 of Goldstein \cite{goldstein} the coupling exists; that
is, there is a probability space
$(\Omega\times\Omega,\mathcal{F}\otimes\mathcal{F},\widetilde
{\mathbb{P}})$
such that
$\widetilde{\mathbb{P}}(\Omega\times\cdot)=\mathbb{P}_x(\cdot)$ and
$\widetilde{\mathbb{P}}(\cdot\times\Omega)=\mathbb{P}_\eta(\cdot)$
(note the dependence on $x$ of $\widetilde{\mathbb{P}}$ is omitted for
notational simplicity). The process on this space is written
$\{\check{X}_n,\check{Y}_n\}_{n\geq0}$ and $T$ is the associated
coupling time. Choose $q \in\mathbb{N}$. For any $\delta>0$,
$M\in\mathbb{N}$, $n>M \vee q$, one has that
\[
\widetilde{\mathbb{P}}\Bigl(\sup_{k \geq n}|S_{k,\check{X}}^{\odot
q}(f) - S_{k,\check{Y}}^{\odot q}(f)|>\delta\Bigr)
\leq
\widetilde{\mathbb{P}}\Bigl(\sup_{k \geq n}|S_{k,\check{X}}^{\odot
q}(f) - S_{k,\check{Y}}^{\odot q}(f)|>\delta, T\leq M\Bigr) +
\widetilde{\mathbb{P}}(T> M)
\]
with $S_{n,\check{X}}^{\odot q}(f)$ as defined in \eqref{eq:ustat_new}.
Now let $\varepsilon>0$ be given and choose $M$ such that
$\widetilde{\mathbb{P}}(T> M)<\varepsilon/2$. The first term on the
right-hand side of the above inequality is now dealt
with:
\begin{eqnarray}\label{eq:sum_prf_neww}
 \hspace{-34pt}&&\widetilde{\mathbb{P}}\Bigl(\sup_{k \!\geq\!
n}|S_{k,\check{X}}^{\odot q}(f) \!-\! S_{k,\check{Y}}^{\odot
q}(f)|\!>\!\delta,
T\!\leq\! M\!\Bigr)
\!=\!\sum_{l=1}^{M}\widetilde{\mathbb{P}}\Bigl(\sup_{k \!\geq\!
n}|S_{k,\check{X}}^{\odot q}(f) \!-\! S_{k,\check{Y}}^{\odot
q}(f)|\!>\!\delta,
T \!=\! l\!\Bigr).
\end{eqnarray}
Then, on the event $\{T=l\}$, one has that the terms involved in the
definitions of $S_{n,\check{X}}^{\odot q}(f)$ and
$S_{n,\check{Y}}^{\odot q}(f)$ only differ for $\vartheta$'s such that
$\vartheta(i)\in\{0,\dots,l-1\}$ for some $i\in\{1,\dots,q\}$. For any
$k>m>0,$ introduce the subset of $\langle q,k+1\rangle$
\[
\Xi_{m,k} := \bigl\{\vartheta\in\langle q,k+1\rangle\dvtx \exists i \in\{1,\ldots,q\}
\mbox{ s.t. } \vartheta(i)< m\bigr\}.
\]
Then for any $l\in\{1,\ldots,M\}$, with $\bar{X}_{\vartheta(i)} =
\check{X}_{\vartheta(i)}$ and $\bar{Y}_{\vartheta(i)} =
\check{Y}_{\vartheta(i)} \mathbb{I}_{\{\vartheta(i)<l\}} +
\check{X}_{\vartheta(i)}\mathbb{I}_{\{\vartheta(i)\geq l\}}$ and the
notation
\[
\Delta(f)_{\bar{X},\bar{Y}}(\vartheta(1),\ldots,\vartheta(q)) :=
f\bigl(\bar{X}_{\vartheta(1)},\dots,\bar{X}_{\vartheta(q)}\bigr) -
f\bigl(\bar{Y}_{\vartheta(1)},\dots,\bar{Y}_{\vartheta(q)}\bigr),
\]
we have
\begin{eqnarray*}
&&
\widetilde{\mathbb{P}}\Bigl(\sup_{k \geq n}|S_{k,\check{X}}^{\odot
q}(f) - S_{k,\check{Y}}^{\odot q}(f)|>\delta, T = l\Bigr)
\\
&&\quad=
\widetilde{\mathbb{P}}\Biggl(\sup_{k \geq
n}\frac{1}{(k+1)_q}\biggl|\sum_{\vartheta\in\Xi_{l,k}}
\Delta(f)_{\bar{X},\bar{Y}}(\vartheta(1),\ldots,\vartheta
(q))\biggr|>\delta,
T = l\Biggr).
\end{eqnarray*}
Let us denote for $l,n \in\mathbb{N}$ such that $n>l$
\[
A_{l,n}:=\biggl\{ \sup_{k \geq n}\frac{1}{(k+1)_q}\biggl|\sum_{\vartheta
\in\Xi_{l,k}}
\Delta(f)_{\bar{X},\bar{Y}}(\vartheta(1),\ldots,\vartheta
(q))\biggr|>\delta
\biggr\} .
\]
It is now shown that $\widetilde{\mathbb{P}}(A_{l,n})$ vanishes as $n
\rightarrow\infty$, which in turn will prove that the above vanishes
as well for any $l\in\{1,\ldots,M\}$. Since $f \in
\mathscr{L}_{(W^r)^{(q)}}$, there exists some (deterministic) constant
$\bar{M}<\infty$ such that
%
\[
\widetilde{\mathbb{P}}(A_{l,n}) \leq\widetilde{\mathbb{P}}\Biggl(
\sup_{k \geq n} \frac{\bar{M}}{(k+1)_q} \Biggl|\sum_{\vartheta\in
\Xi_{l,k}} \Biggl\{\sum_{i=1}^q\bigl[W\bigl(\bar{X}_{\vartheta(i)}\bigr)^r +
W\bigl(\bar{Y}_{\vartheta(i)}\bigr)^r\bigr]\Biggr\}\Biggr|>\delta\Biggr).
\]
Consequently
\begin{eqnarray*}
\widetilde{\mathbb{P}}(A_{l,n})
&\leq&
\mathbb{P}_x\Biggl( \sup_{k \geq n}\frac{\bar{M}}{(k+1)_q} \Biggl|\sum_{\vartheta\in\Xi_{l,k}}
\Biggl\{\sum_{i=1}^q\bigl[W\bigl(X_{\vartheta(i)}\bigr)^r + W\bigl(X_{\vartheta(i)}\bigr)^r\mathbb{I}_{\{\vartheta(i)\geq l\}}\bigr]\Biggr\}\Biggr|>\delta/2\Biggr)
\\
&&{}+
\mathbb{P}_\eta\Biggl( \sup_{k \geq n} \frac{\bar{M}}{(k+1)_q}\Biggl|\sum_{\vartheta\in\Xi_{l,k}} \Biggl\{\sum_{i=1}^q
W\bigl(Y_{\vartheta(i)}\bigr)^r\mathbb{I}_{\{\vartheta(i)<l\}}\Biggr\}\Biggr|>\delta/2\Biggr).
\end{eqnarray*}
The drift condition on $P$ yields the classical result $\sup_{i\geq0}
\{ \mathbb{E}_x[W(X_i)] \vee\mathbb{E}_\eta[W(Y_i)] \} < \infty$. Note
in addition that the cardinality of $\Xi_{l,k}$ is
\[
l\pmatrix{k\cr q-1}q! = (k+1)_q \frac{ql}{k+1} .
\]
Hence one may use an $\mathbb{L}_p$-proof similar to that in
Proposition \ref{poisf1}, with $p\in(1,1/r)$ along with a
Borel--Cantelli argument via Markov's inequality, to conclude that
$\lim_{n\rightarrow\infty} \widetilde{\mathbb{P}}(A_{l,n})=0$. This
allows us to complete the proof by choosing $n$ such that each of the
$M$ terms in the summation \eqref{eq:sum_prf_neww} is less than
$\varepsilon/2M$.

\end{pf}

\begin{sllnu}\label{lemma:vstatistics}
Assume \textup{(A2)} and \textup{(A3)}. Let $q\in\mathbb{N}$, $r\in[0,1)$, $f\in\mathscr{L}_{
(W^r)^{(q)} }$, $x \in E$ and $\{X_i\}$ be as defined earlier. Then,
\[
\lim_{n\rightarrow\infty}|[S_{n,X}^{\otimes q} - \eta^{\otimes q}](f)|
\rightarrow0\qquad \mathbb{P}_x\textrm{-a.s.}
\]
\end{sllnu}

\begin{pf}
The idea of the proof is to use the almost sure convergence results for
$U$-statistics of ergodic stationary processes established in
\cite{aaronson}. In order to achieve this, the coupling~%
$\widetilde{\mathbb{P}}$ introduced in Proposition
\ref{appendixprop:markovcoupl} is utilized. In particular, for any
$\delta>0,$ consider the following upper bound
\begin{eqnarray*}
&&
\mathbb{P}_x\Bigl(\sup_{k\geq n} |[S_{k,X}^{\otimes q} - \eta^{\otimes
q}](f)|>\delta\Bigr)
\\
&&\quad=
\widetilde{\mathbb{P}}\Bigl(\sup_{k\geq n}|[S_{k,\check{X}}^{\otimes q} - S_{k,\check{X}}^{\odot q}](f)
+
[S_{k,\check{X}}^{\odot q} - S_{k,\check{Y}}^{\odot q}](f)
+
[S_{k,\check{Y}}^{\odot q} - \eta^{\otimes q}](f)|>\delta\Bigr)
\\
&&\quad\leq
\mathbb{P}_x\Bigl(\sup_{k\geq n}|[S_{k,X}^{\otimes q} -S_{k,X}^{\odot q}](f)|>\delta/3\Bigr)
+
\widetilde{\mathbb{P}}\Bigl(\sup_{k\geq n}|[S_{k,\check{X}}^{\odot q} -S_{k,\check{Y}}^{\odot q}](f)|>\delta/3\Bigr)
\\
&&\qquad{}+
\mathbb{P}_{\eta}\Bigl(\sup_{k\geq n}|[S_{k,Y}^{\odot q} -\eta^{\otimes q}](f)|>\delta/3\Bigr).
\end{eqnarray*}
%
The convergence to zero of terms on the right-hand side of the
inequality above from right to left are now considered.
Since $\{Y_n\}_{n\geq0}$ is an homogeneous Markov chain, started
in stationarity, it is a stationary ergodic process. In addition, as
$f$ is bounded by integrable products, $(E,\mathcal{E})$ is Polish and
$\{Y_n\}_{n\geq0}$ is absolutely regular (or weakly Bernoulli)
\cite{doukan}, Theorem U of \cite{aaronson} can be invoked; the last
term goes to zero (note that the proofs of \cite{aaronson} extend to
Polish spaces). By Proposition~\ref{appendixprop:markovcoupl}, the
second term goes to zero.

Let us turn to the first term on the right-hand side of the inequality
above. We use an argument similar to that of Theorem 5.1 of
\cite{grams}. This uses the following identity
\[
(n+1)^q [S_{n,X}^{\odot q} - S_{n,X}^{\otimes q}](f) = [(n+1)^q -
(n+1)_{q}]S_{n,X}^{\odot q}(f) -
\sum_{\vartheta\in\overline{\langle q,n+1\rangle}}f\bigl(X_{\vartheta(1)},\ldots
,X_{\vartheta(q)}\bigr),
\]
where $\overline{\langle q,n+1\rangle}:=(q,n+1)\setminus\langle q,n+1\rangle$. Let $p\in
(1,1/r)$. Since $f\in\mathscr{L}_{(W^r)^{(q)}}$, for any
$(i_1,\dots,i_q) \in\{0,\dots,n\}^q$ then by Minkowski's inequality,
followed by Jensen's inequality and the fact that via the drift
condition $\sup_{i \geq0}\mathbb{E}_x[W(pr)(X_i)]<M W(pr)(x)$ for
some $M<\infty$
\[
\mathbb{E}_{x}[|f(X_{i_1},\dots,X_{i_q})|^p]^{1/p}\leq
\|f\|_{(W^r)^{(q)}}\sum_{l=1}^q \mathbb
{E}_x[W(X_{i_l})^{rp}]^{1/p}\leq
M q \|f\|_{W^{(q)}} W^r(x).
\]
As a result
\[
\mathbb{E}_{x}[|S_{n,X}^{\odot q}(f)|^p]^{1/p} \leq M q \|f\|_{W^{(q)}}
W^r(x)
\]
and
\[
\mathbb{E}_{x}\Biggl[\biggl|\sum_{\vartheta\in\overline
{\langle q,n+1\rangle}}f\bigl(X_{\vartheta(1)},\dots,X_{\vartheta(q)}\bigr)\biggr|^p\Biggr]^{1/p}
\leq M [(n+1)^q - (n+1)_q] q \|f\|_{W^{(q)}} W^r(x),
\]
which allows us to conclude that there exists $C_q < \infty$ such that
for any $n > q$
\[
\mathbb{E}_x[(n+1)^q|[S_{n,X}^{\odot q} - S_{n,X}^{\otimes
q}](f)|^p]^{1/p} \leq C_q [(n+1)^q - (n+1)_{q}]W^r(x).
\]
Now since $(n+1)^q - (n+1)_{q} = \mathrm{O}(n^{q-1})$ and $p>1,$ a
Borel--Cantelli argument can be used. The proof of the lemma now
follows.
\end{pf}

\section{Verifying the assumptions}\label{sec:verifyproof}

\begin{pf*}{Proof of Proposition \ref{prop:verify}}
Verifying many of the assumptions
(A1) and (A2) is fairly simple
and can be found in, for example, \cite{jarner}
(i.e.,~(A1)(i)(iii) and (A2)). The small-set condition
(A1)(ii) can easily
be proved in a similar way to the proof of Theorem 2.2 in~%
\cite{roberts2} and is thus omitted. This leaves us with the latter
part of (A2) ((A3) is clearly
true here).\looseness=-1

In our case,
\[
V(x) = \biggl[\frac{|\pi|_{\infty}}{\pi(x)}\biggr]^{s_{v}}
\]
for any $s_w\in(0,1)$ (see \cite{jarner}, Theorems 4.1 and 4.3). The
expression for $W(x)$
\[
\biggl[\frac{|\pi|_{\infty}}{\pi(x)}\biggr]^{\tilde{\alpha}s_{w}},\qquad
s_w\in(0,1),
\]
follows similarly. For the last part of (A2), fix
$r^*,s_w\in(0,1);$ then
\[
\frac{V(x)}{W(x)^{r^*}} = |\pi|_{\infty}^{s_v-r^*\tilde{\alpha}
s_w}\pi(x)^{r^*\tilde{\alpha} s_w-s_v},
\]
which is upper bounded if $s_v\in(0,r^*\tilde{\alpha} s_w)$.
\end{pf*}
\end{appendix}

\subsection*{Acknowledgements}

We thank two referees and the associate editor for valuable comments
that vastly improved the paper. We also thank the editor for his
patience with the paper. The second and third authors acknowledge the
Institute of Statistical Mathematics, Japan, for their support during
the writing of the paper. We also thank Adam Johansen for some useful
comments on previous versions.

\printhistory

\end{document}